\documentclass[10pt]{article}

\usepackage{amsmath,amssymb,latexsym,amsthm,epsfig,stmaryrd}
\usepackage[OT1,T1]{fontenc}

\newtheorem{thm}{Theorem}[section]
\newtheorem{definition}{Definition}[section]
\newtheorem{prp}[thm]{Proposition}
\newtheorem{lem}[thm]{Lemma}
\newtheorem{cor}[thm]{Corollary}

\newenvironment{rmk}{\begin{trivlist}\item[]{\bf Remark:}\setlength{\parindent}{0pt}}
{\end{trivlist}}
\newenvironment{ex}{\begin{trivlist}\item[]{\bf Example:}\setlength{\parindent}{0pt}}
{\end{trivlist}}
\newenvironment{prf}{\begin{trivlist}\item[]{\bf Proof: }}
{\hfill $\blacksquare$ \end{trivlist}}

\def\note#1{\marginpar{\raggedright\if@twoside\ifodd\c@page\raggedleft\fi\fi\sf\scriptsize RMK: #1}}
\def\End{\mathop{\rm End}\nolimits}

\def\ker{\mathop{\rm ker}\nolimits}

\def\Re{\mathop{\rm Re}\nolimits}

\def\deg{\mathop{\rm deg}\nolimits}
\def\rk{\mathop{\rm rk}\nolimits}

\def\cliff{\mbox{\sl Cliff}}

\def\mod{\mathop{\rm mod}\nolimits}

\newcommand{\eunorm}[1]{\parallel\mbox{\hspace{-3pt}} #1\mbox{\hspace{-3pt}}\parallel}

\newcommand{\htimes}{\widehat{\otimes}}
\newcommand{\R}{\mathbb{R}}
\newcommand{\C}{\mathbb{C}}

\newcommand{\mf}{\mathfrak}
\newcommand{\mb}{\mathbf}
\newcommand{\mc}{\mathcal}

\newcommand{\be}[1]{\begin{equation}\label{#1}}
\newcommand{\ee}{\end{equation}}
\newcommand{\gtilde}{\widetilde{\mc{G}}}

\begin{document}
\title{Calibrated cycles and T--duality}
\author{Florian Gmeiner and Frederik Witt}
\date{}
\maketitle

\centerline{\textbf{Abstract}}\noindent
For a subclass of Hitchin's generalised geometries we introduce and analyse the concept of a structured submanifold which encapsulates the classical notion of a calibrated submanifold. Under a suitable integrability condition on the ambient geometry, these generalised calibrated submanifolds minimise a functional occurring as D--brane energy in type II string theories. Further, we investigate the behaviour of calibrated cycles under T--duality and construct non--trivial examples.

\tableofcontents
\section{Introduction}

In this paper, we introduce and investigate a notion of structured submanifold for a subclass of Hitchin's {\em generalised geometries}.  The most basic examples are provided by Harvey's and Lawson's calibrated submanifolds~\cite{hala82}, a notion which has been intensively studied in Riemannian geometry, see for instance~\cite{jo04} and references therein. 

Recall that on a Riemannian manifold $(M,g)$, a differential form $\rho$ of pure degree $p$ is said to define a {\em calibration}, if for any oriented $p$--dimensional subspace $L\subset T_xM$ with induced Riemannian volume form $\varpi_L$, the inequality
\begin{equation}\nonumber
g_x(\rho_x,\varpi_L)\leq\eunorm{\varpi_L}_g=1
\end{equation}
holds ($\eunorm{\cdot}_g$ denoting the norm on $\Lambda^pT^*M$ induced by $g$). Submanifolds whose tangent spaces meet the bound everywhere are said to be {\em calibrated} by $\rho$. Examples comprise K\"ahler manifolds with $\rho$ the K\"ahler form, for which the calibrated submanifolds are precisely the complex submanifolds. More generally, there are natural calibrations for all geometries on Berger's list, and calibrated submanifolds constitute a natural class of structured submanifolds for these. An important property of calibrated submanifolds is to be locally homologically volume--minimising if the calibration is closed. This links into string theory, where the notion of a calibrated submanifold can be adapted to give a geometric interpretation of {\em D--branes} in type I theory, extended objects that extremalise a certain energy functional (cf. for instance~\cite{gip02},~\cite{gupa99},~\cite{mmms00} or the survey~\cite{gmw04}).

On the other hand, geometries defined by forms are the starting point of generalised geometry as introduced in Hitchin's foundational article~\cite{hi03}. The r\^ole of $p$--forms is now played by even and odd forms which we can interpret as {\em spinor fields} for the bundle $TM\oplus T^*\!M$ endowed with its natural orientation and inner product which is contraction. The distinctive feature of this setup is that it can be acted on by both diffeomorphisms and 2--forms or {\em B--fields} $B\in\Omega^2(M)$: These map $X\oplus\xi\in TM\oplus T^*\!M$ to $X\oplus(X\llcorner B/2+\xi)$ and spinor fields $\rho\in\mb{S}_{\pm}\cong\Omega^{ev,od}(M)$ to \begin{equation}\nonumber
e^B\wedge\rho=(1+B+B\wedge B/2+\ldots)\wedge\rho.
\end{equation} 
In particular, this action preserves the natural bilinear form $\langle\cdot\,,\cdot\rangle$ on $\mb{S}=\mb{S}_+\oplus\mb{S}_-\cong\Omega^*(M)$ and transforms closed spinor fields into closed ones if $B$ itself is closed.

Following this scheme, a generalised calibration will be an even or odd form which we view as a $TM\oplus T^*\!M$--spinor field $\rho$. The r\^ole of the Riemannian metric is now played by a {\em generalised Riemannian metric} $\mc{G}$, an involution on $TM\oplus T^*\!M$ compatible with the inner product. Further, a generalised Riemannian metric also induces a norm $\eunorm{\!\cdot\!}_{\mc{G}}$ on $\mb{S}_{\pm}$. Both $\mc{G}$ and the norm transform naturally under $B$--fields. However, it is not altogether clear what the $B$--field action ought to be on a submanifold. The solution to this is, taking string theory as guidance, to consider pairs $(L,F)$ consisting of an embedded submanifold $j_L:L\hookrightarrow M$ and a {\em closed} $2$--form $F\in\Omega^2(L)$. A closed $B$--field $B$ then acts on this by $(L,F)\mapsto(L,F+j^*_LB)$. With such a pair, we can associate the {\em pure} spinor field (cf. Section~\ref{gencalplanes}) 
\begin{equation}\nonumber
\tau_{L,F}=e^F\wedge\theta^{k+1}\wedge\ldots\wedge\theta^n\in\mb{S}_{|L},
\end{equation} 
where the $\theta^i\in\Omega^1(M)_{|L}$ define a frame of the annihilator $N^*\!L\subset T^*\!M$ of $L$. The calibration condition we adopt is
\begin{equation}\nonumber
\langle\rho,\tau_{L,F}\rangle\leq\eunorm{\tau_{L,F}}_{\mc{G}},
\end{equation} 
which therefore makes essential use of the generalised Riemannian metric\footnote{For various other notions of generalised submanifolds in generalised geometry, where no metric structure is involved, see for instance~\cite{bebo04} in the mathematical and~\cite{cgj04} and~\cite{za04} in the physics literature.}.

The first important result we prove is this (cf. Theorem~\ref{spincrit}):
\vspace{-7pt}
\begin{itemize}
	\item Any spinor field of the form $\rho^{ev,od}=e^B\wedge\Re\big([\Psi_L\otimes\Psi_R]^{ev,od}\big)$ defines a calibration, where $[\Psi_L\otimes\Psi_R]$ is the complex form naturally associated with the bi--spinor $\Psi_L\otimes\Psi_R\in\Delta\otimes\Delta$, for $\Delta$ the spinor bundle associated with $TM$.
	\item For such a $\rho$, the pair $(L,F)$ is calibrated if and only if 
	\begin{equation}\nonumber
\mc{A}(\Psi_L)=c\frac{e^{-B}\wedge\tau_{L,F}}{\eunorm{e^{-B}\wedge\tau_{L,F}}_{\mc{G}}}\cdot\Psi_R,
	\end{equation} 
where $\mc{A}$ denotes the so--called {\em charge conjugation operator}, $c$ a complex constant and~$\cdot$ the usual Clifford action of $TM$ on $\Delta$. 
\end{itemize}
This type of calibrations defines what we call a $G_L\times G_R$--{\em structure}. Generalised $G_2$--, $Spin(7)$-- and $SU(m)$--structures~\cite{wi04},~\cite{wi05},~\cite{wi06}, which are natural generalisations of the geometries of holonomy $G_2$, $Spin(7)$ and $SU(m)$ to Hitchin's setup, are examples of these. This result also substantially extends work of Dadok and Harvey~\cite{daha93},~\cite{ha91} on classical calibrations defined by spinors. 

Further, calibrated pairs $(L,F)$ locally minimise the functional
\begin{equation}\nonumber
\mc{E}(L,F)=\int_L\eunorm{\tau_{L,F}}_{\mc{G}}
\end{equation}
if the calibration form is closed. As we will explain, this integral represents the so--called {\em Dirac--Born--Infeld} term of the D--brane energy. Considering an inhomogeneous equation of the form $d\rho=\varphi$ also accounts for the so--called {\em Wess--Zumino} term. We discuss the relationship of calibrated pairs $(L,F)$ with D--branes in {\em type II} string theory as we go along. In particular, our definition of a generalised calibration embraces the cases discussed in~\cite{ko05},~\cite{ma06},~\cite{masm05} in the physics literature. 

The final part of the paper is devoted to T--duality. In string theory, this duality interchanges type IIA with type IIB theory and plays a central r\^ole in the SYZ--formulation of Mirror symmetry~\cite{syz96}. It is also of considerable mathematical interest (cf. for instance~\cite{bem04},~\cite{brs05} and~\cite{busc04}). 
Locally, T--duality is enacted according to the Buscher rules~\cite{bu87} which have a natural implementation in generalised geometry. Our final result (Theorem~\ref{caltdual}) states that if the entire data are invariant under the flow of the vector field along which we T--dualise, the calibration condition is stable under T--duality (as it is expected from a physical viewpoint). We will use this to construct some non--trivial examples of calibrated pairs.

Finally, we remark that our definition of a generalised calibration also makes sense for {\em twisted} generalised geometries~\cite{hi03}, as we shall explain in the paper.

The outline of the paper is this. Based on~\cite{hi03} and~\cite{hi05}, we introduce the setup of generalised geometry in Section~\ref{gengeos}. We thereby considerably extend ideas developed in~\cite{wi06} by generalising the theory to $G_L\times G_R$--structures in arbitrary dimension and by incorporating twisted structures. Section~\ref{gencalcycles} gives the definition of a generalised calibration and calibrated pairs. We investigate their properties and discuss first examples. Section~\ref{physics} makes the link with string theory. Finally, Section~\ref{tduality} is devoted to T--duality and gives further examples.

\begin{center}
{\em Acknowledgments}
\end{center}

\noindent We wish to thank Claus Jeschek for his collaboration in the early stages of this project. The second author was a member of the SFB 647 ``Space.Time.Matter'' funded by the DFG. He also would like to thank Paul Gauduchon and Andrei Moroianu for helpful discussions and creating a very enjoyable ambiance during his stay at \'Ecole Polytechnique, Palaiseau. Further, he thanks the Max--Planck--Institut f\"ur Physik, M\"unchen, for a kind invitation during the preparation of this paper. 
%
%
%
%
%
\section{Generalised geometries}
\label{gengeos}
%
\subsection{The generalised tangent bundle}
%
\subsubsection{Basic setup}\label{setup}
Let $M^n$ be an $n$--dimensional manifold. The underlying geometric setup we will use throughout this paper is known under the name of {\em generalised geometry} and was introduced in Hitchin's foundational paper~\cite{hi03}. Its key feature is to be acted on by so--called $B$--field transformations induced by $\Omega^2(M)$. To implement these, we need to pass from the tangent bundle $TM$ to the bundle $TM\oplus T^*\!M$. We endow $TM\oplus T^*\!M$ with (a) its canonical orientation and (b) the metric of signature $(n,n)$ given by contraction, namely\footnote{Our conventions slightly differ from~\cite{hi03} and~\cite{wi04} which results in different signs and scaling factors.} 
\begin{equation}\label{innerprod}
(X\oplus\xi,X\oplus\xi)=X\llcorner\xi=\xi(X).
\end{equation} 
Therefore, this vector bundle is associated with a principal $SO(n,n)$--fibre bundle $P_{SO(n,n)}$. As $GL(n)\subset SO(n,n)$, $P_{SO(n,n)}$ is actually obtained as an extension of the frame bundle $P_{GL(n)}$ associated with $TM$. Now take $B\in\Omega^2(M)$ and think of it as a linear map $B:TM\to T^*\!M$ sending $X$ to the contraction $X\llcorner B=B(X,\cdot)$. Then we define the corresponding $B$--field transformation by
\begin{equation}\nonumber
X\oplus\xi\in TM\oplus T^*\!M\mapsto X\oplus\big(\xi+\frac{1}{2}X\llcorner B\big).
\end{equation}
This transformation can be interpreted as an exponential: For a point $q\in M$, $\big(T_qM\oplus T^*_q\!M,(\cdot\,,\cdot)_q\big)$ is an oriented pseudo--Euclidean vector space of signature $(n,n)$. Its symmetry group is $SO(n,n)$, whose Lie algebra we may identify with
\begin{equation}\nonumber
\mf{so}(n,n)\cong\Lambda^2(T_qM\oplus T^*_q\!M)\cong\Lambda^2T_qM\oplus T_q^*\!M\otimes T_qM\oplus\Lambda^2T^*_q\!M.
\end{equation}
Hence, $B_q\in\Lambda^2T^*_q\!M$ becomes a skew--symmetric endomorphism of $T_qM\oplus T^*_q\!M$ under the identification $\zeta\wedge \eta (X\oplus\xi)=(\zeta,X\oplus\xi) \xi-(\eta,X\oplus\xi)\zeta=X\llcorner(\zeta\wedge\eta)/2$. Applying the exponential map $\exp:\mf{so}(n,n)\to SO(n,n)$, we obtain
\begin{equation}\nonumber
\exp(B_q)(X\oplus\xi)=(\sum\limits_{l=0}^{\infty}\frac{B_q^l}{l!})(X\oplus\xi)=X\oplus\xi+\frac{1}{2}X\llcorner B_q
\end{equation}
as $B^2_q(X\oplus\xi)=B_q(X\llcorner B_q)/2=0$. With respect to the splitting $T_qM\oplus T^*_q\!M$, we will also use the matrix representation
\begin{equation}\nonumber
e^{B_q}=\left(\begin{array}{cc} Id & 0\\ \frac{1}{2}B_q & Id\end{array}\right).
\end{equation}

Further, $TM\oplus T^*\!M$ admits a natural bracket which extends the usual vector field bracket $[\cdot\,,\cdot]$, the so--called {\em Courant bracket}~\cite{co90}. It is defined by
\begin{equation}\nonumber
\llbracket X\oplus\xi,Y\oplus\eta\rrbracket=[X,Y]+\mc{L}_X\eta-\mc{L}_Y\xi-\frac{1}{2}d(X\llcorner\eta-Y\llcorner\xi),
\end{equation}
where $\mc{L}$ denotes the Lie derivative. If $B$ is a {\em closed} $2$--form, then $\exp(B)$ commutes with $\llbracket\cdot\,,\cdot\rrbracket$.

The inclusion $GL(n)\leqslant SO(n,n)$ can be lifted, albeit in a non--canonical way, to $Spin(n,n)$. It follows that the $SO(n,n)$--structure is always spinnable. In the sequel, we assume all base manifolds to be orientable, hence the canonic lift of the inclusion\footnote{The notation $G_+$ refers to the identity component of a given Lie group $G$.} $GL(n)_+\hookrightarrow SO(n,n)$
induces a canonic spin structure $P_{GL(n)_+}\subset P_{Spin(n,n)_+}$ with associated spinor bundle
\begin{equation}\nonumber
\mb{S}(TM\oplus T^*\!M)=P_{GL(n)_+}\times_{GL(n)_+}S.
\end{equation} 
In order to construct the spin representations $S$ for $Spin(n,n)$, we consider its vector representation $\R^{n,n}$ which we split into two maximal isotropic subspaces, $\R^{n,n}=W\oplus W'$. Using the inner product $g=g^{n,n}$, we may identify $W'$ with $W^*$. Under this identification, $g(w',w)=w'(w)/2$, that is, we recover precisely~(\ref{innerprod}). Now $X\oplus\xi\in W\oplus W^*$ acts on $\rho\in\Lambda^*W^*$ via
\begin{equation}\nonumber
(X\oplus\xi)\bullet\rho=-X\llcorner\rho+\xi\wedge\rho.
\end{equation}
This action squares to minus the identity, and by the universal property of Clifford algebras, it extends to an algebra isomorphism 
\begin{equation}\label{cliffrep}
\cliff(n,n)\cong\cliff(W\oplus W^*)\cong\End(\Lambda^*W^*).
\end{equation}
Our initial choice of a splitting into isotropic subspaces is clearly immaterial, as $\cliff(n,n)$ is a {\em simple} algebra, hence the representation~(\ref{cliffrep}) is unique up to isomorphism. The spin representation of $Spin(n,n)$ is therefore isomorphic with $S=\Lambda^*W^*$ and can be invariantly decomposed into the modules of {\em chiral} spinors 
\begin{equation}\nonumber
S_{\pm}=\Lambda^{ev,od}W^*.
\end{equation} 
Moreover, after choosing some trivialisation of $\Lambda^nW^*$, we can define the bilinear form 
\begin{equation}\nonumber
\langle\rho,\tau\rangle=[\rho\wedge\widehat{\tau}]^n\in\Lambda^nW^*\cong\R.
\end{equation}
Here, $[\,\cdot\,]^n$ denotes projection on the top degree component and $\widehat{\,\cdot\,}$ is the sign--changing operator defined on forms of degree $p$ by 
\begin{equation}\nonumber
\widehat{\alpha^p}=(-1)^{p(p+1)/2}\alpha^p.
\end{equation} 
Then
\begin{equation}\nonumber
\langle (X\oplus\xi)\bullet\rho,\tau\rangle=(-1)^n\langle\rho,(X\oplus\xi)\bullet\tau\rangle
\end{equation} 
and in particular, this form is $Spin(n,n)_+$--invariant. It is symmetric for $n\equiv 0,3\mod4$ and skew for $n\equiv 1,2\mod4$, i.e. 
\begin{equation}\nonumber
\langle\rho,\tau\rangle=(-1)^{n(n+1)/2}\langle\tau,\rho\rangle.
\end{equation}
Moreover, $S^+$ and $S^-$ are non--degenerate and orthogonal if $n$ is even and totally isotropic if $n$ is odd.

We also have an induced action of the Lie algebra $\mf{so}(n,n)$. By exponentiation, $B\in\Lambda^2W^*$ acts on $\rho\in S_{\pm}$ via
\begin{equation}\nonumber
e^B\bullet\rho=(1+B+\frac{1}{2}B\wedge B+\ldots)\wedge\rho.
\end{equation}
This exponential links into the $B$--field transformation on $W\oplus W^*$ via the 2-1 covering map $\pi_0:Spin(n,n)\to SO(n,n)$, namely
\begin{equation}\label{coveringB}
\pi_0(e^B_{Spin(n,n)})=e^{\pi_{0*}(B)}_{SO(n,n)}=e^{2B}_{SO(n,n)}=\left(\begin{array}{cc}Id & 0\\B & Id\end{array}\right).
\end{equation}
On the other hand, the element $\sum A^m_l\xi^l\otimes X_m\in\mf{gl}(n)\subset\mf{so}(n,n)$ acts on $\rho$ via 
\begin{equation}\nonumber
\sum A^m_l\xi^l\otimes X_m(\rho)=\frac{1}{2}\sum A^m_l[\xi^l,X_m]\bullet\rho=\frac{1}{2}{\rm Tr}(A)+A^*\rho,
\end{equation} 
where $A^*\rho$ denotes the natural extension to forms of the dual representation of $\mf{gl}(n)$. Consequently, 
\begin{equation}\nonumber
S_{\pm}\cong\Lambda^{ev,od}W^*\otimes\sqrt{\Lambda^nW}\mbox{ as a $GL(n)_+$--space,}
\end{equation}
and we obtain for the spinor bundle
\begin{equation}\nonumber
\begin{array}{ccccl}\mb{S}_{\pm} & = & \mb{S}(TM\oplus T^*\!M)_{\pm} & \cong & P_{GL(n)_+}\times_{GL(n)_+}\Lambda^{ev,od}\R^{n*}\otimes\sqrt{\Lambda^n\R^n}\\[5pt]
 & & & = & \Lambda^{ev,od}T^*\!M\otimes\sqrt{\Lambda^nTM}.
\end{array}
\end{equation}
The choice of a trivialisation of the line bundle $\sqrt{\Lambda^nTM}$, that is, of a nowhere vanishing $n$--vector field $\nu$, induces thus an isomorphism between $TM\oplus T^*\!M$--spinor fields and even or odd differential forms. Put differently, it reduces the structure group of $M$ from $GL(n)_+$ to $SL(n)$, for which the bundles $\mb{S}_{\pm}$ and $\Lambda^{ev,od}T^*\!M$ are isomorphic. We shall denote this isomorphism by $\mf{L}^{\nu}$ to remind ourselves that it depends on the choice of the $n$--vector $\nu$.
Moreover, from a $GL(n)_+$--point of view, the bilinear form $\langle\cdot\,,\cdot\rangle$ takes values in the reals, as for $\rho,\,\tau\in S\cong\Lambda^*W^*\otimes\sqrt{\Lambda^nW}$, we have $\langle\rho,\tau\rangle\in\Lambda^nW^*\otimes\sqrt{\Lambda^nW}\otimes\sqrt{\Lambda^nW}\cong\R$.
In particular, using the isomorphism $\mf{L}^{\nu}$, we get
\begin{equation}\nonumber
\langle\rho,\tau\rangle=\nu\big([\mf{L}^{\nu}(\rho)\wedge\widehat{\mf{L}^{\nu}(\tau)}]^n\big)\in C^{\infty}(M).
\end{equation}
\subsubsection{Twisting with an $H$--flux}
\label{hfluxtwist}
More generally still, there is a twisted version of the previously discussed setup which is locally modeled on $TM\oplus T^*\!M$, leading to the notion of a {\em generalised tangent bundle}~\cite{hi05}. 

Let $H\in\Omega^3(M)$ be a closed $3$--form (in physicists' terminology, this is the so--called $H$--{\em flux}). Choose a convex cover of coordinate neighbourhoods $\{U_a\}$ of $M$, whose induced transition functions are $s_{ab}:U_{ab}=U_a\cap U_b\to GL(n)_+$. Locally, 
\begin{equation}\nonumber
H_{|U_a}=dB^{(a)}
\end{equation} 
for some $B^{(a)}\in\Omega^2(U_a)$. Then we can define the $2$--forms
\begin{equation}\nonumber \beta^{(ab)}=B^{(a)}_{|U_{ab}}-B^{(b)}_{|U_{ab}}\in\Omega^2(U_{ab})
\end{equation} 
which are {\em closed} (this will be of importance later). On $U_{ab}$, we can use the coordinates provided by either the chart $U_a$ or $U_b$, and consider $\beta^{(ab)}$ as a map from $U_{ab}\to\Lambda^2\R^{n*}$; we denote this map by $\beta^{(ab)}_a$ or $\beta^{(ab)}_b$. Any two trivialisations are therefore related by $\beta^{(ab)}_a=s_{ab}^*\beta^{(ab)}_b$ and similarly for all other trivialised sections. Out of the transition functions
\begin{equation}\nonumber
S_{ab}=\left(\begin{array}{cc} s_{ab} & 0\\ 0 & s_{ab}^{-1{\rm tr}}\end{array}\right)\in GL(n)_+\subset SO(n,n)_+
\end{equation}
of the bundle $TM\oplus T^*\!M$, we define new transition functions
\begin{equation}\nonumber
\sigma_{ab}=S_{ab}\circ e^{2\beta^{(ab)}_b}=e^{2\beta^{(ab)}_a}\circ S_{ab}:U_{ab}\to SO(n,n)_+.
\end{equation}
These satisfy indeed the {\em cocycle condition} on $U_{abc}\not=\emptyset$, i.e.
\begin{eqnarray}
\sigma_{ab}\circ \sigma_{bc} & = & S_{ab}\circ e^{2\beta^{(ab)}_b}\circ e^{2\beta^{(bc)}_b}\circ S_{bc}\nonumber\\
& = & S_{ab}\circ e^{2(B^{(a)}_b-B^{(b)}_b+B^{(b)}_b-B^{(c)}_b)}\circ S_{bc}\nonumber\\
& = & S_{ab}\circ S_{bc}\circ e^{2(B^{(a)}_c-B^{(c)}_c)}\nonumber\\
& = & S_{ac}\circ e^{2\beta^{(ac)}_c}\nonumber\\
& = & \sigma_{ca}^{-1}.\label{cocycle}
\end{eqnarray}
Therefore, we can define a vector bundle, the {\em generalised tangent bundle}, by
\begin{equation}\nonumber
\mb{E}=\mb{E}(H)=\prod U_a\times(\R^n\oplus\R^{n*})/\sim_{\sigma_{ab}},
\end{equation}
where two triples $(a,p,X\oplus\xi)$ and $(b,q,Y\oplus\eta)$ are equivalent if and only if $p=q$ and $X\oplus\xi=\sigma_{ab}(p)(Y\oplus\xi)$. A section $t\in\Gamma(\mb{E})$ can therefore be regarded as a family of smooth maps $\{t_a:U_a\to\R^n\oplus\R^{n*}\}$ transforming under $s_a=\sigma_{ab}(t_b)$. As the notation suggests, the generalised tangent bundle depends, up to isomorphism, only on the closed $3$--form $H$. Indeed, assume we are given a different convex cover $\{U_a'\}$ together with locally defined $2$--forms $B^{(a)'}\in\Omega^2(U'_a)$ such that $H_{|U'_a}=dB^{(a)'}$. This results in a new family of transition functions $\sigma'_{ab}=S_{ab}\circ\exp(2\beta^{(ab)'}_b)$. Now on the intersection $V_a=U_a\cap U'_a$, the locally defined $2$--forms $\widetilde{G}^{(a)}=B^{(a)}_{|V_a}-B^{(a)'}_{|V_a}$ are closed, and one readily verifies the family $G_a=\exp(\widetilde{G}^{(a)}_a)$ to define a gauge transformation, i.e.
\begin{equation}\nonumber
\sigma'_{ab}=G^{-1}_a\circ \sigma_{ab}\circ G_b\mbox{ on }V_{ab}\not=\emptyset.
\end{equation}
In particular, the bundles defined by the families $\sigma_{ab}$ and $\sigma'_{ab}$ respectively are isomorphic. The sections of the bundle $\mb{E}(H)$ relate to $\Gamma(TM\oplus T^*\!M)$ via the map $\mf{K}:\Gamma(TM\oplus T^*M)\to\Gamma\big(\mb{E}(H)\big)$ defined by
\begin{equation}\nonumber
\mf{K}_a(X_a\oplus\xi_a)=e^{2B^{(a)}_a}(X_a\oplus\xi_a)=X_a\oplus B^{(a)}_a(X_a)+\xi_a.
\end{equation}
This is indeed well--defined as
\begin{eqnarray*}
\sigma_{ba}\mf{K}_a(X_a\oplus\xi_a) & = & e^{\beta^{(ba)}_b}\!\big(X_b\oplus B^{(a)}_b(X_b)+\xi_b\big)\\
& = & X_b\oplus\big(B^{(b)}_b(X_b)+\xi_b\big)\\
& = & \mf{K}_b\big(S_{ba}(X_a\oplus\xi_a)\big).
\end{eqnarray*}
Since the transition functions $\sigma_{ab}$ take values in $SO(n,n)$, the invariant orientation and inner product $(\cdot\,,\cdot)$ on $\R^n\oplus\R^{n*}$ render $\mb{E}$ an oriented pseudo--Riemannian vector bundle. Furthermore, twisting with the closed $2$--forms $\beta^{(ab)}$ preserves the Courant bracket, which is therefore also defined on $\mb{E}(H)$. Since any other choice of a trivialisation leads to a gauge transformation induced by {\em closed} $B$--fields, this additional structure does not depend on the choice of a trivialisation either. For the pseudo--Riemannian structure on $\mb{E}$ thus defined, we can also consider spinor fields. A canonic spin structure is given by \begin{equation}\nonumber
\widetilde{\sigma}_{ab}=\widetilde{S}_{ab}\bullet e^{\beta^{(ab)}_b}=e^{\beta^{(ab)}_a}\bullet\widetilde{S}_{ab},
\end{equation} 
where $\widetilde{S}_{ab}\in GL(n)_+\subset Spin(n,n)_+$ denotes the lift of $S_{ab}$ and where we exponentiate $\beta^{(ab)}$ to $Spin(n,n)_+$, so that $\pi_0\circ\widetilde{\sigma}_{ab}=\sigma_{ab}$. The even and odd spinor bundles associated with $\mb{E}$ are
\begin{equation}\nonumber
\mb{S}(\mb{E})_{\pm}=\coprod\limits_aU_a\times S_{\pm}/\sim_{\widetilde{\sigma}_{ab}}.
\end{equation}
An $\mb{E}$--spinor field $\rho$ is thus represented by a collection of smooth maps $\rho_a:U_a\to S_{\pm}$ with $\rho_a=\widetilde{\sigma}_{ab}\bullet\rho_b$. 
In order to make contact with differential forms, let $\nu$ be a trivialisation of $\Lambda^nTM$ which we think of as a family of maps $\nu_a=\lambda_a^{-2}\nu_0\in\Lambda^n\R^n$ with $\lambda_a^{-2}\in C^{\infty}(M)$. It follows that $\lambda_b=\sqrt{s_{ab}}\cdot\lambda_b$. The map $\mf{L}^{\nu}_{\mb{E}}:\Gamma\big(\mb{S}(\mb{E})_{\pm}\big)\to\Omega^{ev,od}(M)$ induced by
\begin{equation}\nonumber
\mf{L}^{\nu}_{\mb{E}a}:(\rho_a:U_a\to S_{\pm})\mapsto (e^{-B^{(a)}_a}\wedge\lambda_a\cdot\rho_a:U_a\to\Lambda^{ev,od}\R^{n*})
\end{equation}
is an isomorphism, and we usually drop the subscript $\mb{E}$ to ease notation. This transforms correctly under the action of the transition functions $s_{ab}$ on $\Lambda^{ev,od}T^*\!M$, as one can show by using the fact that $\rho_a=\widetilde{\sigma}_{ab}\bullet\rho_b$. Indeed, over $U_{ab}$ we have
\begin{eqnarray*}
s_{ab}^*(e^{-B^{(b)}_b}\wedge\lambda_b\cdot\rho_b) & = & \lambda_a\sqrt{\det s_{ab}}\cdot e^{-B^{(b)}_a}\wedge s^*_{ab}\rho_b\nonumber\\
& = & \lambda_a\cdot e^{-B^{(a)}_a}\wedge e^{\beta^{(ab)}_a}\bullet \widetilde{S}_{ab}\bullet\rho_b\nonumber\\
& = & \lambda_a\cdot e^{-B^{(a)}_a}\wedge\rho_a,
\end{eqnarray*}
or equivalently, $\mf{L}^{\nu}_a\circ\widetilde{\sigma}_{ab}=s_{ab}^*\circ \mf{L}^{\nu}_b$. The operators $\mf{L}^{\nu}$ and $\mf{K}$ relate via
\begin{equation}\nonumber
\mf{L}^{\nu}\big(\mf{K}(X\oplus\xi)\bullet\rho\big)=-X\llcorner\mf{L}^{\nu}(\rho)+\xi\wedge\mf{L}^{\nu}(\rho)=:(X\oplus\xi)\bullet\mf{L}^{\nu}(\rho).
\end{equation}
In the same vein, the $Spin(n,n)$--invariant form $\langle\cdot\,,\cdot\rangle$ induces as above a globally defined inner product on $\Gamma(\mb{S})$ by $\langle\rho,\tau\rangle=\nu\big([\mf{L}^{\nu}(\rho)\wedge\widehat{\mf{L}^{\nu}(\tau)}]^n\big)$.

Twisting with closed $B$--fields allows the definition of a further operator, namely
\begin{equation}\nonumber
d_{\nu}:\Gamma\big(\mb{S}(E)_{\pm}\big)\to\Gamma\big(\mb{S}(E)_{\mp}\big),\quad(d_{\nu}\rho)_a=\lambda^{-1}_a\cdot d_a(\lambda_a\cdot\rho_a),
\end{equation}
where $d_a$ is the usual differential applied to forms $U_a\to\Lambda^*\R^{n*}$, i.e. $(d\alpha)_a=d_a\alpha_a$. This definition gives indeed rise to an $\mb{E}$--spinor field, for
\begin{eqnarray*}
\widetilde{\sigma}_{ab}\bullet(d_{\nu}\rho)_b & = & \lambda^{-1}_b\cdot e^{\beta^{(ab)}_a}\wedge \sqrt{\det s_{ab}}\cdot s_{ab}^*d_b(\lambda_b\cdot\rho_b)\\
& = & \lambda_a^{-1}\cdot d_a\big(e^{\beta^{(ab)}_a}\wedge s_{ab}^*(\lambda_b\cdot\rho_b)\big)\\
& = & (d_{\nu}\rho)_a.
\end{eqnarray*}
The operator $d_{\nu}$ squares to zero and therefore induces an elliptic complex on $\Gamma\big(\mb{S}(E)_{\pm}\big)$. It computes the so--called {\em twisted cohomology}, where on replaces the usual differential $d$ of de Rham cohomology by the twisted differential $d_H=d+H\wedge$.

\begin{prp}\label{twisteddiff}
Let $\rho\in\Gamma\big(\mb{S}(E)\big)$. Then
\begin{equation}\nonumber
\mf{L}^{\nu}(d_{\nu}\rho)=d_H\mf{L}^{\nu}(\rho).
\end{equation}
\end{prp}

\noindent The assertion follows from a straightforward local computation.

\begin{rmk}
The identification of $\mb{S}(E)_{\pm}$ with even or odd forms depends on the choice of the local $B$--fields $B^{(a)}$. However, any other choice leads, as we have seen, to a gauge transformation by closed $2$--forms $\widetilde{G}^{(a)}$, which therefore preserves the closeness of the induced differential forms.
\end{rmk}
%
\subsection{Generalised metrics}
\label{genmetrics}
In addition to the $3$--form flux $H$ and a nowhere vanishing $n$--vector field $\nu\in\Gamma(\Lambda^nTM)$, we now wish to build in a Riemannian metric to define {\em generalised Riemannian metrics}. 

Again, we first consider the bundle $TM\oplus T^*\!M$. Let us think of a Riemannian metric as a map $g:TM\to T^*\!M$, which is invertible in virtue of the non--degeneracy. With respect to the decomposition $TM\oplus T^*\!M$, we define the endomorphism
\begin{equation}\nonumber
\mc{G}_0=\left(\begin{array}{cc} 0 & g^{-1}\\
g & 0\end{array}\right).
\end{equation}
The $B$--transform of $\mc{G}_0$ induced by $B\in\Omega^2(M)$ is
\begin{eqnarray*}
\mc{G}_B & = & e^{2B}\circ\mc{G}_0\circ e^{-2B}\\
& = & \left(\begin{array}{cc} Id & 0\\
B & Id\end{array}\right)\circ\left(\begin{array}{cc} 0 & g^{-1}\\
g & 0\end{array}\right)\circ\left(\begin{array}{cc} Id & 0\\
-B & Id\end{array}\right)\\
& = & \left(\begin{array}{cc} -g^{-1}B & g^{-1}\\ g-Bg^{-1}B & Bg^{-1}\end{array}\right).
\end{eqnarray*}
To simplify notation, we will usually write $\mc{G}$ (unless we want to emphasise the $B$--field) and refer to $\mc{G}$ as the {\em generalised metric} induced by $g$ and $B$. Note that $\mc{G}$ squares to the identity, and its $\pm1$--eigenspaces $V^{\pm}$ give a {\em metric splitting}, i.e. an orthogonal decomposition
\begin{equation}\nonumber 
TM\oplus T^*\!M=V^+\oplus V^-.
\end{equation}
into maximally positive/negative definite subbundles $V^{\pm}$. The restriction of $(\cdot\,,\cdot)$ to $V^{\pm}$ will be denoted by $g_{\pm}$. Conversely, any metric splitting gives rise to an honest Riemannian metric $g$ and $B\in\Omega^2(M)$, whose corresponding generalised metric has $V_{\pm}$ as $\pm1$--eigenspaces~\cite{wi04}. In terms of structure groups, we note that the global decomposition of $TM\oplus T^*\!M$ into $V^+\oplus V^-$ gives rise to a reduction from $SO(n,n)_+$ to
\begin{eqnarray*}
e^{2B}SO(n,0)\times SO(0,n)e^{-2B} & \cong & \{\left(\begin{array}{cc} A_+ & 0\\ 0 & A_-\end{array}\right)\,|\,A_{\pm}:V^{\pm}\to V^{\pm}\mbox{ such that}\\
& &\phantom{\{}A^*_{\pm}g_{\pm}=g_{\pm}\mbox{ and }\det A_{\pm}=1\}\\
 & = & SO(V^+)\times SO(V^-),
\end{eqnarray*}
where $SO(n,0)\times SO(0,n)$ is the subgroup preserving the eigenspace decomposition of $\mc{G}_0$.

\begin{definition}{\rm (\cite{hi05},~\cite{wi04})} A {\em generalised (Riemannian) metric} for the generalised tangent bundle $\mb{E}(H)$ is the choice of a maximally positive definite subbundle $V^+$.
\end{definition}

\noindent Put differently, a generalised metric induces a splitting of the exact sequence
\begin{equation}\label{sequence}
0\to T^*\!M\to\mb{E}\to TM\to 0.
\end{equation}
We denote the lift of vector fields $X\in\Gamma(TM)$ to sections in $\Gamma(V^+)$ by $X^+$. Locally, $X^+$ corresponds to smooth maps $X_a^+:U_a\to\R^n\oplus\R^{n*}$ with $X^+_a=\sigma_{ab}(X^+_b)$, and $X^+_a=X_a\oplus P^+_aX_a$ for linear isomorphisms $P^+_a:\R^n\to\R^{n*}$. From the transformation rule on $\{X^+_a\}$, we deduce 
\begin{equation}\nonumber
\beta^{(ab)}_a=P_a^+-s_{ab}^*P^+_bs_{ba}
\end{equation}
(where $X_a=s_{ab}(X_b)$). As above, the symmetric part $g_a=(P_a^++P_a^{+{\rm tr}})/2$ is positive definite, and since $\beta^{(ab)}$ is skew--symmetric, the symmetrisation of the right hand side vanishes. Hence
\begin{equation}\nonumber
\frac{1}{2}(P_a^+-s_{ab}^*P^+_bs_{ba}+P_a^{+{\rm tr}}-s_{ab}^*P^{+{\rm tr}}_bs_{ba})=g_a-s_{ab}^*g_bs_{ba}=0,
\end{equation}
so that the collection $g_a:U_a\to\odot^2\R^{n*}$ of positive definite symmetric $2$--tensors patches together to a globally defined metric. Conversely, a Riemannian metric $g$ induces a generalised Riemannian structure on $\mb{E}(H)$: The maps $P_a=B^{(a)}_a+g_a$ induce local lifts of $TM$ to $\mb{E}$ which give rise to a global splitting of~(\ref{sequence}).

\begin{prp}
A generalised Riemannian structure is characterised by the datum $(g,H)$, where $g$ is an honest Riemannian metric and $H$ a closed $3$--form. 
\end{prp}

\begin{rmk}
Of course, the negative definite subbundle $V^-$ also defines a splitting of~(\ref{sequence}). The lift of a vector field $X$ is then induced by $X^-_a=X_a\oplus P^-_aX_a$ with $P^-_a=-g_a+B^{(a)}_a$.
\end{rmk}

\noindent In presence of a metric, we can pick the canonic $n$--vector field locally given by $\nu_{g|U_a}=(\det g_{|U_a})^{-1/2}\partial_{x^1}\wedge\ldots \wedge\partial_{x^n}$. We shall write $\mf{L}$ for $\mf{L}^{\nu_g}$. From our original choice $\nu$ this differs by a scalar function, i.e. $\nu=e^{2\phi}\nu_g$ for $\phi\in C^{\infty}(M)$, so that $\mf{L}^{\nu}=e^{-\phi}\mf{L}$. We then write $\mf{L}^{\phi}=\mf{L}^{\nu}$ and $d_{\phi}=d_{\nu}$. For reasons becoming apparent later, $\phi$ is referred to as the {\em dilaton field}.

A generalised metric also induces further structure on spinor fields. On every fibre $T_qM\oplus T^*_q\!M$, the straight generalised metric $\mc{G}_0$ acts as a composition of reflections. Namely, if $e_1,\ldots,e_n$ is an orthonormal basis of $T_qM$, then $V^{\pm}_q$ is spanned by $v_k^{\pm}=e_k\oplus\pm g(e_k)$ and
\begin{equation}\nonumber
\mc{G}_0=R_{v^-_1}\circ\ldots\circ R_{v^-_n}\in O(n,n),
\end{equation}
where $R_{v^-_k}$ denotes reflection along $v_k^-$. Hence its lift $\gtilde_0$ to $Pin(n,n)$ is given, up to a sign, by the Riemannian volume form
\begin{equation}\nonumber \gtilde_0=\varpi_{V^-_p}=v^-_1\bullet\ldots\bullet v^-_n
\end{equation}
on $V^-_q$. The operator $\gtilde_0$ is therefore globally well--defined. Note that $\mc{G}_0$ preserves or reverses the natural orientation on $TM\oplus T^*\!M$ if $n$ is even or odd, so that $\gtilde_0$ preserves or reserves the chirality of $TM\oplus T^*\!M$--spinor fields accordingly. There is an alternative description of its action which will turn out to be useful. For this, we denote by $\mf{J}$ the natural isomorphism between $\cliff(TM,g)$ and $\Lambda^*T^*\!M$. Recall that for any $X\in TM$ and $a\in \cliff(TM,g)$ of pure degree,
\begin{equation}\nonumber
\mf{J}(X\cdot a)=-X\llcorner\mf{J}(a)+X\wedge\mf{J}(a),\quad\mf{J}(a\cdot X)=(-1)^{\deg(a)}\big(X\llcorner\mf{J}(a)+X\wedge\mf{J}(a)\big).
\end{equation}
Moreover,
\begin{equation}\nonumber
\star_g\mf{J}(a)=\mf{J}(\widehat{a}\cdot\varpi_g)
\end{equation}
where $\varpi_g$ denotes the Riemannian volume form on $TM$ as well as its image in $\cliff(TM,g)$ under $\mf{J}$ by abuse of notation. If $\sim$ is the involution defined by $\widetilde{a}^{ev,od}=\pm a^{ev,od}$, then $\varpi_g\cdot a=a\cdot\varpi_g$ if $n$ is odd, while $\varpi_g\cdot a=\pm \widetilde{a}\cdot\varpi_g$ for $n$ even. By computing 
\begin{equation}\nonumber
\mf{L}\big(d_1^-\bullet\ldots\bullet d_n^-\bullet\mf{J}(\mf{J}^{-1}(\rho))\big)=\pm\mf{L}\big(\mf{J}(\omega_g\cdot\mf{J}^{-1}(\rho))\big),
\end{equation}
where the sign depends on $n=\rk TM$, we deduce that $\mf{L}(\gtilde_0\rho)$ is equal to $\star_g\widehat{\mf{L}(\rho)}$ for $n$ even and $\pm\star_g\widehat{\mf{L}(\rho)}$ for $n$ odd and $\rho$ of even or odd parity. For a $B$--field transformed Riemannian metric with $V^{\pm}=\exp(2B)(D^{\pm})$, we conjugate $\mc{G}_0$ by $\exp(2B)$ and therefore $\varpi_{D^-}$ by $\exp(B)$ in view of~(\ref{coveringB}). Hence $\varpi_{V^-}\bullet\rho=e^B\bullet \varpi_{D^-}\bullet e^{-B}\bullet\rho$. Up to signs, $\gtilde$ therefore coincides with the $\Box$--operator in~\cite{wi04}. The twisted case follows easily, for instance if $n$ is even,
\begin{eqnarray*}\nonumber
\mf{L}_a(\varpi_{V^-a}\bullet\rho_a) & = & \mf{L}_a(e^{B^{(a)}_a}\bullet \varpi_{D^-a}\bullet e^{-B^{(a)}_a}\bullet\rho_a)\\
& = & \varpi_{D^-a}\bullet\mf{L}_a(\rho_a)\\
& = & \star_{g_a}\widehat{\mf{L}_a(\rho_a)}.
\end{eqnarray*}

\begin{prp}
If $V^+$ defines a generalised Riemannian metric, then the action of $\widetilde{\mc{G}}=\varpi_{V^-}$ on $\mb{S}(\mb{E})_{\pm}$ is given by
\begin{equation}\nonumber
\mf{L}(\gtilde\rho)=\left\{\begin{array}{ll}$n$\mbox{ even:}\quad&\star_g\widehat{\mf{L}(\rho)}\\$n$\mbox{ odd:}\quad&\star_g\widehat{\widetilde{\mf{L}(\rho)}}\end{array}\right.=:\gtilde\big(\mf{L}(\rho)\big),
\end{equation}
where $g$ is the Riemannian metric induced by $V^+$. 
\end{prp}

\noindent Note that 
\begin{equation}\nonumber
\gtilde^2=(-1)^{n(n+1)/2}Id\quad\mbox{ and}\quad\langle\gtilde\rho,\tau\rangle=(-1)^{n(n+1)/2}\langle\rho,\gtilde\tau\rangle.
\end{equation} 
In particular, $\widetilde{\mc{G}}$ is an isometry for $\langle\cdot\,,\cdot\rangle$ and defines a complex structure on $\mb{S}(\mb{E})$ if $n\equiv1,2\mod4$. Moreover, given a generalised metric, we define the bilinear form 
\begin{equation}\nonumber
\mc{Q}^{\pm}(\rho,\tau)=\pm(-1)^m\langle\rho,\gtilde\tau\rangle\mbox{ on }\mb{S}(\mb{E})_{\pm},\quad n=2m,2m+1.
\end{equation}
Then for $n=2m$ (and similarly for $n=2m+1$),
\begin{eqnarray*}
\mc{Q}^{\pm}(\rho,\tau) & = & \pm(-1)^m\langle\rho,\gtilde\tau\rangle\\& = & \pm(-1)^m\nu_g[\mf{L}(\rho)\wedge\big(\star_g\widehat{\mf{L}(\tau)}\big)^{\wedge}]^n\\ 
& = & g\big(\mf{L}(\rho),\mf{L}(\tau)\big),
\end{eqnarray*}
where we used the general rule $\star\widehat{\alpha}^{ev,od}=\pm(-1)^m\widehat{\star\alpha}^{ev,od}$. We shall denote the associated norm on $\mb{S}_{\pm}$ by $\eunorm{\cdot}_{\mc{G}}$.

If the manifold is spinnable, the presence of a generalised metric also implies a very useful description of the complexified spinor modules $\mb{S}(\mb{E})_{\pm}\otimes\C$ as the tensor product $\Delta_n(TM)\otimes\Delta_n(TM)$ of the complex spin representation $\Delta_n$ of $Spin(n,0)$. From standard representation theory, we have algebra isomorphisms
\begin{equation}\label{spiniso}
\cliff^{\,\C}(TM,\pm g)\cong\cliff(TM\otimes\C,g^{\C})=\left\{\begin{array}{lc}$n$\mbox{ even:}\quad&\End_{\C}(\Delta_n)\\
$n$\mbox{ odd:}\quad&\End_{\C}(\Delta_n)\oplus\End_{\C}(\Delta_n)\end{array}\right..
\end{equation} 
The module $\Delta_n$ is the space of {\em (Dirac) spinors} and, as a complex vector space, is isomorphic with $\C^{2^{[n/2]}}$. Further, $\Delta_n$ carries an hermitian inner product $q$ for which $q(a\cdot\Psi,\Phi)=q(\Psi,\widehat{a}\cdot\Phi)$, $a\in\cliff(TM\otimes\C,g^{\C})$. By convention, we take the first argument to be conjugate--linear. Restricting the isomorphism~(\ref{spiniso}) to the Spin groups of signature $(p,q)$ yields the complex spin representation $Spin(p,q)\to GL(\Delta_{p+q})$. If $p+q$ is even, then there is a decomposition $\Delta_{p+q}=\Delta_{p+q,+}\oplus\Delta_{p+q,-}$ into the irreducible $Spin(p,q)$--representations $\Delta_{p+q,\pm}$, the so--called {\em Weyl spinors} of {\em positive} and {\em negative chirality}. Finally, in all dimensions, there exists a conjugate--linear endomorphism $\mc{A}$ of $\Delta_n$ (the {\em charge conjugation operator} in physicists' language) such that~\cite{wa89}
\begin{equation}\nonumber 
\mc{A}(X\cdot\Psi)=(-1)^{m+1}X\cdot \mc{A}(\Psi)\mbox{ and }\mc{A}^2=(-1)^{m(m+1)/2}Id,\quad n=2m,2m+1.
\end{equation} 
In particular, $\mc{A}$ is $Spin(n)$--equivariant. Moreover, $\mc{A}$ reverses the chirality for $n=2m$, $m$ odd. 
We define a $Spin(n)$--invariant bilinear form (which, abusing of notation, we also write $\mc{A}$) by
\begin{equation}\nonumber
\mc{A}(\Psi,\Phi)=q(\mc{A}(\Psi),\Phi),
\end{equation} 
for which
\begin{equation}\nonumber
\mc{A}(\Psi,\Phi)=(-1)^{m(m+1)/2}\mc{A}(\Phi,\Psi)\mbox{ and }\mc{A}(X\cdot\Psi,\Phi)=(-1)^m\mc{A}(\Psi,X\cdot\Phi)
\end{equation}
if $n=2m,\,2m+1$. We can inject the tensor product
$\Delta_n\otimes\Delta_n$ into $\Lambda^*\C^n$ 
by sending $\Psi_L\otimes\Psi_R$ to the form of mixed degree 
\begin{equation}\nonumber
[\Psi_L\otimes\Psi_R](X_1,\ldots,X_n)=\mc{A}\big(\Psi_L,(X_1\wedge\ldots\wedge X_n)\cdot\Psi_R\big).
\end{equation} 
In fact, this is an isomorphism for $n$ even. In the odd case, we obtain an isomorphism by concatenating $[\cdot\,,\cdot]$ with projection on even or odd forms, which we write as $[\cdot\,,\cdot]^{ev,od}$. Since this map is $Spin(n)$--equivariant, it acquires global meaning over $M$, and we use the same symbol for the resulting map $\Delta_n(TM)\otimes\Delta_n(TM)\to\Omega^{ev,od}(M)$ (referred to as the {\em fierzing} map in the physics' literature). Next we define
\begin{equation}\nonumber
[\cdot,\cdot]^{\mc{G}\,ev,od}:\Gamma\big(\Delta_n(TM)\otimes\Delta_n(TM)\big)\stackrel{[\cdot,\cdot\,]^{ev,od}}{\longrightarrow}\Omega^{ev,od}(M)\otimes\C\stackrel{\mf{L}^{-1}}{\longrightarrow}\Gamma\big(\mb{S}(\mb{E})_{\pm}\otimes\C\big).
\end{equation}
A vector field $X$ acts on $TM$--spinor fields via the inclusion $TM\hookrightarrow\cliff(TM,g)$ and Clifford multiplication. On the other hand, we can lift $X$ to sections $X^{\pm}$ of $V^{\pm}$ which act on $\mb{E}$--spinor fields via the inclusion $V^{\pm}\hookrightarrow\cliff(TM\oplus T^*\!M)$ and Clifford multiplication. To see how these actions are related by the map $[\cdot\,,\cdot]^{\mc{G}}$, we use the fact that the orthogonal decomposition of $TM\oplus T^*\!M$ into $V^+\oplus V^-$ makes $\cliff(TM\oplus T^*\!M)$ isomorphic with the twisted tensor product $\cliff(V^+)\htimes\cliff(V^-)$\footnote{The twisted tensor product $\htimes$ of two graded algebras $\mf{A}$ and $\mf{B}$ is defined on elements of pure degree as $a\htimes b\cdot a'\htimes b'=(-1)^{\deg(b)\cdot\deg(a')} a\cdot a'\htimes b\cdot b'$}. The vector bundles $(V^{\pm},g_{\pm})$ are isometric to $(TM,\pm g)$ via the lifts $X\in\Gamma(TM)\mapsto X^{\pm}\in\Gamma(V^{\pm})$. By extending $X\htimes Y\mapsto X^+\bullet Y^-$ to $\cliff(TM\oplus T^*\!M)$, we get a further isomorphism 
\begin{equation}\nonumber
\cliff^{\,\C}(TM,g)\htimes\cliff^{\,\C}(TM,-g)\cong\cliff^{\,\C}(V^+\oplus V^-,g_+\oplus g_-).
\end{equation} 

\begin{prp}\label{commutprop}
We have
\begin{eqnarray*}\nonumber
\,[X\cdot\Psi_L\otimes\Psi_R]^{\mc{G}} & = & (-1)^{n(n-1)/2}X^+\bullet[\Psi_L\otimes\Psi_R]^{\mc{G}},\\
\,[\Psi_L\otimes Y\cdot\Psi_R]^{\mc{G}} & = & -Y^-\bullet\widetilde{[\Psi_L\otimes\Psi_R]^{\mc{G}}}.
\end{eqnarray*}
\end{prp}

\begin{prf}
Fix a point $q\in U_a$ and an orthonormal basis $e_1,\ldots,e_n$ of $T_qM$. To ease notation, we drop any reference to $q$ and $U_a$ and assume first $B_a=0$. By definition and the properties of $\mc{A}$ recalled above,
\begin{eqnarray*}
[e_j\cdot\Psi_L\otimes\Psi_R]  & = & \sum\limits_{K} q\big(\mc{A}(e_j\cdot\Psi_L),e_K\cdot\Psi_R\big)e_K\\
& = & (-1)^{n(n-1)/2+1}\sum\limits_{K}q\big(\mc{A}(\Psi_L),e_j\cdot e_K\cdot\Psi_R\big)e_K\\
& = & (-1)^{n(n-1)/2+1}\sum\limits_{K}q\big(\mc{A}(\Psi_L),(-e_j\llcorner e_K+e_j\wedge e_K)\cdot\Psi_R\big)e_K\\
& = & (-1)^{n(n-1)/2}\big(\sum\limits_{j\in K}q\big(\mc{A}(\Psi_L),e_j\llcorner e_K\cdot\Psi_R\big)e_j\wedge (e_j\llcorner e_K)-\\
& & \phantom{(-1)^{n(n-1)/2}\big(}\sum\limits_{j\not\in K}q\big(\mc{A}(\Psi_L),e_j\wedge e_K\cdot\Psi_R\big)e_j\llcorner (e_j\wedge e_K)\big)\\
& = & (-1)^{n(n-1)/2}\big(-e_j\llcorner+g(e_j)\wedge\big)\bullet[\Psi_L\otimes\Psi_R].
\end{eqnarray*} 
Compounding with $\mf{L}^{-1}$, we therefore get
\begin{equation}\nonumber
[X\cdot\Psi_L\otimes\Psi_R]^{\mc{G}_0}=X^+\bullet[\Psi_L\otimes\Psi_R]^{\mc{G}_0}.
\end{equation}
For non--trivial local B--fields $B^{(a)}$, we have 
\begin{eqnarray*}
\,[X\cdot\Psi_L\otimes\Psi_R]^{\mc{G}}_{|U_a} & = & e^{B_a}\bullet\big(X\oplus g(X)\big)\bullet e^{-B_a}\bullet\mf{L}^{-1}[\Psi_L\otimes\Psi_R]_{|U_a}\\
& = & \pi_0(e^{B_a})\big(X\oplus g(X)\big)\bullet [\Psi_L\otimes\Psi_R]^{\mc{G}}_{|U_a}\\
& = & X^+\bullet [\Psi_L\otimes\Psi_R]^{\mc{G}}_{|U_a}.
\end{eqnarray*}
The second assertion follows in the same fashion.
\end{prf}

\noindent A first, but important observation is that $\mb{E}$--spinor fields corresponding to bi--spinor fields are ``self--dual'' in the following sense: The element $\varpi_{V^-}$ in $\cliff(V^+\oplus V^-)$ corresponds to $1\htimes \varpi_g$, so that 
\begin{equation}\nonumber
\gtilde[\Psi_L\otimes\Psi_R]^{\mc{G}}=(-1)^m[\Psi_L\otimes \varpi_g\cdot\Psi_R]^{\mc{G}}
\end{equation} 
for $n=2m$ and 
\begin{equation}\nonumber
\gtilde[\Psi_L\otimes\Psi_R]^{\mc{G}}=(-1)^m\widetilde{[\Psi_L\otimes \varpi_g\cdot\Psi_R]^{\mc{G}}}
\end{equation}
for $n=2m+1$. By standard Clifford representation theory, the action of the Riemannian volume form $\varpi_g$ on chiral spinor fields is given by
\begin{equation}\label{volaction}
\begin{array}{lcrl}
\varpi_g\cdot\Psi_{\pm} & = & \pm (-1)^{m(m+1)/2}i^m\Psi_{\pm},\quad&\mbox{ if }n=2m,\\[5pt] \varpi_g\cdot\Psi & = & (-1)^{m(m+1)/2}i^{m+1}\Psi,\quad&\mbox{ if }n=2m+1.
\end{array}
\end{equation}
From this we deduce the

\begin{cor}\label{selfdual}
Let $\Psi_{L,R}\in\Delta_n$. 

\noindent{\rm (i)} If $n=2m$ and $\Psi_{L,R}$ are chiral, then for $\Psi_R\in\Delta_{\pm}$
\begin{equation}\nonumber
\gtilde[\Psi_L\otimes\Psi_R]^{\mc{G}}=\pm (-1)^{m(m-1)/2}i^m[\Psi_L\otimes\Psi_R]^{\mc{G}}.
\end{equation}
\noindent{\rm (ii)} If $n=2m+1$, then
\begin{equation}\nonumber
\gtilde[\Psi_L\otimes\Psi_R]^{\mc{G}}=(-1)^{m(m-1)/2}i^{m+1}\widetilde{[\Psi_L\otimes\Psi_R]^{\mc{G}}}.
\end{equation}
\end{cor}
%
\subsection{$G_L\times G_R$--structures}
\label{gengstruc}
Within the generalised setting, further reductions can be envisaged. For the theory of generalised calibrations we are aiming to develop in Section~\ref{gencalcycles}, a particular class of generalised structures is of interest. These generalise special geometries defined by $TM$--spinor fields:

\begin{ex}
A $G_2$--{\em structure} on a seven--fold $M^7$ is defined by a $3$--form $\varphi$ with the property that for all $p\in M$, $\varphi_p\in\Lambda^3T_p^*\!M$ lies in the unique orbit diffeomorphic to $GL(7)/G_2$. The structure group therefore reduces from $GL(7)$ to $G_2$. Since $G_2\subset SO(7)$, $\varphi$ induces a Riemannian metric. Conversely, fix a Riemannian metric $g$. If the manifold is spinnable, we can consider the spinor bundle $\Delta(TM)$ associated with $Spin(7)$. Its spin representation $\Delta$ is of {\em real} type, that is, it is the complexification of an irreducible real spin representation $\Delta_{\R}\cong \R^8$. The unit sphere is diffeomorphic to $Spin(7)/G_2$ and therefore, a unit spinor field also defines a $G_2$--structure. This links into the ``form'' definition via the fierzing map $[\cdot\,,\cdot]$, namely
\begin{equation}\label{straightg2}
[\Psi\otimes\Psi]=1-\varphi-\star_g\varphi+\varpi_g.
\end{equation}
\end{ex}

To make contact with generalised geometries, we note that $x\in\C^n\mapsto ix\in\C^n$ extends to a Clifford algebra isomorphism $\cliff(\R^n,g)\otimes\C\cong\cliff(\R^n,-g)\otimes\C$. Restricted to $Spin(n,0)\subset\cliff(\R^n,g)\otimes\C$, this gives  an isomorphism $a\in Spin(n,0)\mapsto\widehat{a}\in Spin(0,n)$. Then $\underline{Y}=Y_1\cdot\ldots\cdot Y_{2l}\in Spin(0,n)$ acts on $\Psi\in\Delta_n$ by $\widehat{\underline{Y}}\cdot\Psi$, so that by Proposition~\ref{commutprop}
\begin{equation}\nonumber
[\Psi_L\otimes\widehat{\underline{Y}}\cdot\Psi_R]^{\mc{G}}=(-1)^l\widehat{\underline{Y}}^-\bullet[\Psi_L\otimes\Psi_R]^{\mc{G}}=\underline{Y}^-\bullet[\Psi_L\otimes\Psi_R]^{\mc{G}}.
\end{equation}

\begin{cor}
$[\cdot\,,\cdot]^{\mc{G}}$ is $Spin(n,0)\times Spin(0,n)$--equivariant.
\end{cor}

\noindent Consider now a Riemannian spin manifold $(M,g)$ and two subgroups $G_{L,R}\subset Spin(n,0)$ which stabilise the collection of spinor fields $\{\Psi_{L,j}\}$ and $\{\Psi_{R,k}\}$ respectively. We identify $G_R$ with its image $\widehat{G}_R$ in $Spin(0,n)$. The tensor product $\Delta(TM)\otimes\Delta(TM)$ is associated with a $Spin(n,0)\times Spin(0,n)$--structure, and the collection of bi--spinor fields $\{\Psi_{L,j}\otimes\Psi_{R,k}\}$ induces a reduction to $G_L\times G_R$ (or more accurately to $G_L\times\widehat{G}_R$, but we will usually omit this) inside $Spin(n,0)\times Spin(0,n)$. By the previous corollary, $[\cdot\,,\cdot]^{\gtilde}$ maps this bi--spinor field to a $G_L\times G_R\subset Spin(V^+)\times Spin(V^-)$--invariant $\mb{E}$--spinor field, where $Spin(V^+)\times Spin(V^-)$ comes from the generalised metric on $\mb{E}$. For instance, taking up the previous example leads to the notion of 

\paragraph{{\em Generalised $G_2$--structures}~\cite{wi04}.}
A {\em generalised $G_2$--structure} can be defined in either of the following three ways:
\vspace{-5pt}
\begin{itemize}
	\item by the datum $(g,H,\Psi_L,\Psi_R)$, where $\Psi_L$, $\Psi_R$ are two real unit spinor fields which reduce the $Spin(7)$--structure $P_{Spin(7)}$ to the principal fibre bundles $P_{G_{2L}}$ and $P_{G_{2R}}$ respectively.
	\vspace{-5pt}
	\item by a principal $G_{2L}\times G_{2R}$--fibre bundle to which the $Spin(7,7)_+$--principal fibre bundle associated with $\mb{E}$ reduces.
	\vspace{-5pt}
	\item by an even or odd $\mb{E}$--spinor fields $\rho$ whose stabiliser is invariant under a group conjugated to $G_{2L}\times G_{2R}$. It is implicitly defined by
	\begin{equation}\nonumber		\mf{L}(\rho)=[\Psi_L\otimes\Psi_R]^{ev,od}.
	\end{equation}
\end{itemize}
The even and odd spinor field are related by $\gtilde(\rho^{ev})=\rho^{od}$
\vspace{-10pt}

\paragraph{{\em Generalised $SU(3)$--structures}~\cite{jewi05a},~\cite{jewi05b}.}
In dimension $6$, we have an isomorphism $Spin(6)\cong SU(4)$ under which the complex spin representations $\Delta_{\pm}$ are isomorphic with the standard vector representations of $SU(4)$ and its complex conjugate, namely $\C^4$ and $\overline{\C^4}$. The unit spheres are isomorphic with $SU(4)/SU(3)$, so that $SU(3)$ stabilises a spinor $\Psi_{\pm}$ in both $\Delta_+$ and $\Delta_-$, which are related by $\Psi_-=\mc{A}(\Psi_+)$. A {\em generalised $SU(3)$--structure} is characterised by either of the following:
\vspace{-5pt}
\begin{itemize}
	\item by the datum $(g,H,\Psi_L,\Psi_R)$, where $\Psi_L$, $\Psi_R$ are two unit spinor fields which reduce the $Spin(6)$--structure $P_{Spin(6)}$ to the principal fibre bundles $P_{SU(3)_L}$ and $P_{SU(3)_R}$ respectively.
	\vspace{-5pt}
	\item by a principal $SU(3)_L\times SU(3)_R$--fibre bundle to which the $Spin(6,6)_+$--principal fibre bundle associated with $\mb{E}$ reduces.
	\vspace{-5pt}
	\item by a pair $(\rho_0,\rho_1)$ of $\mb{E}$--spinor fields whose stabiliser is invariant under a group conjugated to $SU(3)_L\times SU(3)_R$. The spinor fields are implicitly defined by
	\begin{equation}\nonumber	\mf{L}(\rho_0)=[\mc{A}(\Psi_L)\otimes\Psi_R],\quad\mf{L}(\rho_1)=[\Psi_L\otimes\Psi_R].
	\end{equation}
\end{itemize}
\noindent Similarly, one can define generalised $SU(m)$--structures (see also~\cite{wi06}). Applying the operator $\mc{A}\otimes\mc{A}$ to $\Delta\otimes\Delta$ yields the other $SU(3)_L\times SU(3)_R$--invariant pair $\big([\Psi\otimes\mc{A}(\Psi)],[\mc{A}(\Psi)\otimes\mc{A}(\Psi)]\big)$. In the particular case (subsequently referred to as ``straight'') where $\Psi_L=\Psi=\Psi_R$, we obtain
\begin{equation}\nonumber
\mf{L}(\rho_0)=[\mc{A}(\Psi)\otimes\Psi]=e^{-i\omega},\quad\mf{L}(\rho_1)=[\Psi\otimes\Psi]=\Omega^{3,0}
\end{equation}
where $\omega$ is the K\"ahler form and $\Omega^{3,0}$ the invariant $(3,0)$--form of the underlying single $SU(3)$--structure.

\paragraph{{\em Generalised $Spin(7)$--structures}~\cite{wi04}.}
In analogy with dimension $7$, the two chiral $Spin(8)$--representations $\Delta_{\pm}$ are of real type, and the stabiliser of a real chiral unit spinor field is isomorphic with $Spin(7)$. Hence, a {\em generalised $Spin(7)$--structure } is characterised by either of the following:
\vspace{-5pt}
\begin{itemize}
	\item by the datum $(g,H,\Psi_L,\Psi_R)$, where $\Psi_L$, $\Psi_R$ are two real unit spinor fields which reduce the $Spin(8)$--structure $P_{Spin(8)}$ to the principal fibre bundles $P_{Spin(7)_L}$ and $P_{Spin(7)_R}$ respectively.
	\vspace{-5pt}
	\item by a principal $Spin(7)_L\times Spin(7)_R$--fibre bundle to which the $Spin(8,8)_+$--principal fibre bundle associated with $\mb{E}$ reduces.
	\vspace{-5pt}
	\item by an $\mb{E}$--spinor field $\rho$ whose stabiliser is invariant under a group conjugated to $Spin(7)_L\times Spin(7)_R$. It is implicitly defined by
	\vspace{-5pt}
	\begin{equation}\nonumber		\mf{L}(\rho)=[\Psi_L\otimes\Psi_R].
	\end{equation}
\end{itemize}
\vspace{-5pt}
The spinor field $\rho$ is even or odd if the spinors are of equal or opposite chirality, reflecting the fact that there are two conjugacy classes of $Spin(7)$ in $Spin(8)$, each of which stabilises a unit spinor in $\Delta_{\R+}$ or $\Delta_{\R-}$ respectively. The straight case induces an even spinor given by
\begin{equation}\nonumber
\mf{L}(\rho)=[\Psi\otimes\Psi]=1-\Omega+\varpi_g, 
\end{equation}
where $\Omega$ is the $\star_g$--selfdual $Spin(7)$--invariant $4$--form.  
If $\rho$ is odd, the principal $Spin(7)_{\!L}$-- and $Spin(7)_{\!R}$--bundle intersect in a principal $G_2$--bundle. In this case, a generalised $Spin(7)$--structure is locally the $B$--field transform of a $G_2$--structure on $M^8$.
%
%
%
%
%
\section{Generalised calibrations}
\label{gencalcycles}
Next we define a notion of structured submanifold for $G_L\times G_R$--structures as introduced above. In the case of a straight structure, we will recover the so--called calibrated submanifolds which were first considered in Harvey's and Lawson's seminal work~\cite{hala82}. We recall some elements of their theory first. Throughout the Sections~\ref{classcal} and~\ref{gencalplanes}, $T$ will denote a real, oriented, $n$--dimensional vector space.
%
\subsection{Classical calibrations}
\label{classcal}
Let $(T,g)$ be a Euclidean vector space and $\rho\in\Lambda^kT^*$. One says that $\rho$ defines a {\em calibration} if for any oriented $k$--dimensional subspace $L\subset T$ (subsequently also referred to as a $k$--{\em plane}), with induced Euclidean volume form $\varpi_L$, the inequality
\begin{equation}\label{classical}
g(\rho,\varpi_L)\leq\sqrt{g(\varpi_L,\varpi_L)}=1
\end{equation}
holds with equality for at least one $k$--plane. Such a plane is then said to be {\em calibrated} by $\rho$. We immediately conclude from the definition that $\rho$ defines a calibration and calibrates $L$ if and only if $\star_g\rho$ defines a calibration and calibrates $L^{\perp}$. Moreover, let $A\in GL(n)$. The calibration condition is $GL(n)$--equivariant in the following sense: $\rho$ defines a calibration with respect to $g$ if and only if $A^*\rho$ does so with respect to $A^*g$. If $L$ is calibrated for $\rho$, then so is $A(L)$ for $A^*\rho$. In particular, if $\rho$ is $G$--invariant, then so is the calibration condition. Therefore, calibrated subplanes live in special $G$--orbits of the Grassmannian $Gr_k(T)$.

A spinorial approach to calibrations was developed by Dadok and Harvey~\cite{daha93},~\cite{ha91} for $G_2$, $Spin(7)$ and further geometries in dimension $n=8m$. Consider the case of a $G_2$-- and $Spin(7)$--structure on $T$ defined by the unit spinor $\Psi$. Then
\begin{itemize}
	\vspace{-5pt}
	\item any homogeneous component $[\Psi\otimes\Psi]^p$ of $[\Psi\otimes\Psi]\in\Lambda^*T^*$ defines a calibration form
	\vspace{-5pt}
	\item a $k$--plane $L$ is calibrated with respect to $[\Psi\otimes\Psi]^k$ if and only if $\Psi=\varpi_L\cdot\Psi$.
\end{itemize}

If $(M,g)$ is a Riemannian manifold, then a $k$--form $\rho\in\Omega^k(M)$ is called a {\em calibration} if and only if $\rho_p$ defines a calibration on $T_pM$ for any $p\in M$. A $k$--dimensional submanifold $j=j_L:L\hookrightarrow M$ is said to be {\em calibrated}, if $T_pj(L)$ is calibrated for every $p\in j(L)$. If the calibration form is closed, then calibrated submanifolds are locally homologically volume--minimising. Here are examples of interest to us.

\begin{ex}\hfill\newline
\noindent(i) {\em $SU(3)$--structures:}
The powers $\omega^l/l!$ of the K\"ahler form and $\Re\Omega$, the real part of the complex volume form, define calibrations. The submanifolds calibrated by the former are {\em complex} (of real dimension $2l$) and {\em special Lagrangian} for the latter (of dimension $3$). Similar remarks apply to $SU(m)$--structures for arbitrary $m$.

\noindent(ii) {\em $G_2$--structures:}
The $G_2$--invariant $3$--form $\varphi$ defines a calibration form whose calibrated submanifolds are called {\em associative}. The 4--dimensional submanifolds calibrated by $\star\varphi$ are referred to as {\em coassociative}.

\noindent(iii) {\em $Spin(7)$--structures:}
The $Spin(7)$--invariant $4$--form $\Omega$ defines a calibration form. The calibrated submanifolds are the so--called {\em Cayley} submanifolds. 
\end{ex}
%
\subsection{Generalised calibrated planes}
\label{gencalplanes}
To carry this notion over to the generalised setup, one is naturally led to replace $k$--forms by even or odd $T\oplus T^*$--spinors. But what ought to be the analogue of $k$--planes? A natural object on which $B$--fields act are the maximally isotropic subplanes of $TM\oplus T^*\!M$. These can be equivariantly identified with lines of {\em pure} spinors. By definition, purity means that the totally isotropic space
\begin{equation}\nonumber
W_{\tau}=\{x\oplus\xi\in T\oplus T^*\,|\,(x\oplus\xi)\bullet\tau=0\}
\end{equation}
is of maximal dimension, that is $\dim W_{\tau}=n$. The following result is classical.

\begin{prp}{\rm (\cite{ch96} Prop.~III.1.9, p.140)} A spinor $\tau\in S_{\pm}$ is pure if and only if it can be written in the form
\begin{equation}\nonumber
\tau=c\cdot e^F\bullet\varpi,
\end{equation}
where $c\in\R_{\not=0}$, $F\in\Lambda^2T^*$ and $\varpi=\theta^1\wedge\ldots\wedge\theta^l$ with $\theta^i\in T^*$ is a decomposable form in $\Lambda^lT^*$. Further, if $W_{\tau}$ denotes the isotropic subspace in $T\oplus T^*$ corresponding to the line in $S_{\pm}$ spanned by $\tau$, then $\theta^1,\ldots,\theta^l$ is a basis of $W_{\tau}\cap T^*$.
\end{prp}

\begin{definition}\label{rankdef}
The number $\rk\tau=n-\dim(W_{\tau}\cap T^*)$ is called the {\em rank} of $\tau$.
\end{definition}

\begin{rmk}
The rank of a pure spinor is invariant under $B$--field transformations. Indeed
\begin{equation}\nonumber
\dim(W_{\tau}\cap T^*)=\dim\big(e^{2B}(W_{\tau}\cap T^*)\big)=\dim(W_{e^B\bullet\tau}\cap T^*),
\end{equation}
for $\exp(2B)$ acts trivially on $T^*$.
\end{rmk}

\noindent In other words, we can write $W_{\tau}\cap T^*=N^*\!L$, where $N^*\!L$ denotes the annihilator of 
\begin{equation}\nonumber
L={\rm Ann}\,\varpi=\{X\in T\,|\,X\llcorner\varpi=0\}.
\end{equation}
Therefore, the decomposable part of a pure spinor is unique up to a scalar. However, the $2$--form $F$ is only determined up to a $2$--form $F'$ with $F'\wedge\varpi=0$. Put differently, $F'$ belongs to $\mc{N}^2_L=\mc{J}(N^*\!L)\cap\Lambda^2T^*$, the subspace of $2$--forms in the $\Lambda^*T^*$--ideal $\mc{J}(N^*\!L)$ generated by $N^*\!L$. Using the metric $g$, we can identify the orthogonal complement $\mc{N}_L^{2\perp}$ with the image of $\Lambda^2L^*$ under the pull--back map induced by the orthogonal projection $p_L:T\to L$. Hence, modulo rescaling any pure spinor can be uniquely written as $\tau=\exp(F_0)\bullet\varpi$, where $F_0=p_L^*(F)\in p^*_L(\Lambda^2L^*)$ for $L={\rm Ann}\,\varpi$.

\begin{definition}
An {\em isotropic pair} $(L,F)$ consists of a $k$--plane $L\subset T$ and a $2$--form $F\in\Lambda^2L^*$.
\end{definition}

\noindent To render the correspondence between pure spinors and isotropic pairs explicit, we make the following choices: Let again $\varpi_L$ denote the induced Euclidean volume element as well as its pull--back to $T$. We then define the pure spinor $\tau_L=\widehat{\star\varpi_L}$, where the Hodge dual is taken with respect to $T$. The corresponding maximally isotropic subspace is $W_L=L\oplus N^*\!L$. By equivariance, the pure spinor
\begin{equation}\nonumber
\tau_{L,F}=e^{p^*_L(F)}\bullet\widehat{\star\varpi_L}\in S_{(-1)^k},
\end{equation}
then corresponds to $W_{L,F}=\exp(2F)W_L$. To ease notation, we shall simply write $F$ instead of $p^*_L(F)$. For later reference, we let 
\begin{equation}\nonumber
[\tau_{L,F}]=\R_{>0}\tau_{L,F}
\end{equation} 
denote the half--line inside $S_{(-1)^k}$ spanned by $\tau_{L,F}$. We say that the isotropic pair $(L,F)$ is {\em associated} with the pure spinor $\tau$, if $[\tau_{L,F}]=[\tau]$.

\begin{cor}\label{purespinchar}
There is a $Spin(n,n)_+$--equivariant 1-1--correspondence between
\begin{itemize}
	\item isotropic pairs $(L,F)$.
	\item half--lines of spinors $[\tau_{L,F}]$.
	\item maximally isotropic oriented subspaces \begin{equation}\label{isorep}
	W_{L,F}=e^{2F}(L\oplus N^*\!L).
\end{equation}
\end{itemize}
\end{cor}

The definition of a calibration is this.

\begin{definition}\label{caldef}
Let $(g,B)$ define a generalised metric on $T\oplus T^*$. A chiral spinor $\rho\in S_{\pm}$ is called a {\em generalised calibration} if for any isotropic pair $(L,F)$, 
\begin{itemize}
	\item the inequality 
\begin{equation}\label{formcal}
\langle\rho,\tau\rangle\leq\eunorm{\tau}_{\mc{G}}
\end{equation}
holds for one (and thus for all) $\tau\in[\tau_{L,F}]$ and
	\item there exists at least one isotropic pair $(L,F)$ for which one has equality.
\end{itemize}
\noindent Isotropic pairs $(L,F)$ which meet the bound are said to be {\em calibrated} by $\rho$.
\end{definition}

\noindent We usually drop the adjective ``generalised'' and simply speak of a calibration. Condition~(\ref{formcal}) is clearly $Spin(n,n)_+$--equivariant. Hence $\rho$ defines a calibration which calibrates $(L,F)$ if and only if $A\bullet\rho$ defines a calibration which calibrates $(L_A,F_A)$ associated with $A\bullet\tau_{L,F}$. Furthermore, as $\gtilde$ acts as an isometry for $\mc{Q}_{\mc{G}}$, $\rho$ defines a calibration if and only if $\gtilde\rho$ does. If $(L,F)$ is calibrated for $\rho$, then the isotropic pair $(L_{\gtilde},F_{\gtilde})$ associated with $(-1)^{n(n+1)/2}[\gtilde\tau_{L,F}]$ is calibrated for $\gtilde\rho$.

Our first proposition states this calibration condition to be the formal analogue of~(\ref{classical}) and makes the appearance of the generalised metric $(g,B)$ explicit. To harmonise with notation used in Section~\ref{calsub}, we let $j=j_L:L\to T$ denote the injection of the subspace $L$ into $T$. The pull--back $j^*=j^*_L:\Lambda^pT^*\to\Lambda^pL^*$ is then simply restriction to $L$.

\begin{prp}\label{voldef}
A spinor $\rho$ defines a calibration if and only if for any isotropic pair $(L,F)$ with $k=\dim L$, the inequality \begin{equation}\nonumber
[e^{-F}\wedge j^*\rho]^k\leq \sqrt{\det\big(j^*(g+B)-F\big)}\varpi_L
\end{equation}
or equivalently, 
\begin{equation}\nonumber
g(e^{-F}\wedge \rho,\varpi_L)\leq \sqrt{\det\big(j^*(g+B)-F\big)},
\end{equation}
holds\footnote{In the latter inequality, the determinant is computed from the matrix of the bilinear form $j^*(g+B)-F$ with respect to some orthonormal basis of $L$}. The bound is met if and only if $(L,F)$ is calibrated.
\end{prp}

\begin{prf}
Regarding $\rho$ as an exterior form, it follows from the remarks above that
\begin{eqnarray*}
\,[e^{-F}\wedge j^*\rho]^k & = & g(e^{-F}\wedge\rho,\varpi_L)\varpi_L\\
& = & \star[e^{-F}\wedge\rho\wedge\star \varpi_L]^n\varpi_L\\
& = & \star[\rho\wedge\widehat{(e^{F}\wedge\widehat{\star\varpi_L})}]^n\varpi_L\\
& = &
\langle\rho,\tau_{L,F}\rangle\varpi_L.
\end{eqnarray*}
The result is now an immediate consequence of the technical lemma below which implies $\eunorm{\tau_{L,F}}_{\mc{G}}\varpi_L=\sqrt{\det\big(j^*(g+B)-F\big)}\varpi_L$. 
\end{prf}

\begin{lem}\label{det}
Let $L$ be an oriented subspace and $\alpha\in\Lambda^2T^*$. Then
\begin{equation}\nonumber
\sqrt{g(e^{\alpha}\wedge\star\varpi_L,e^{\alpha}\wedge\star\varpi_L)}\varpi_L=\sqrt{\det\big(j^*(g+\alpha)\big)}\varpi_L =\sqrt{\det\big(j^*(g-\alpha)\big)}\varpi_L.
\end{equation}
\end{lem}

\begin{prf}
Fix an oriented orthonormal basis $e_1,\ldots,e_k$ of $L$ such that  $j^*\alpha=\sum_{l=1}^{[k/2]}a_le_{2l-1}\wedge e_{2l}$. Then
\begin{eqnarray}
\sqrt{\det\big(j^*(g-\alpha)\big)} & = & \prod^{[k/2]}_{l=1}\sqrt{1+a_l^2}\label{determinant}\\
& = & \sqrt{1+\sum_{l=1}^{[k/2]}a_l^2+\sum_{l_1<l_2}a_{l_1}^2\cdot a_{l_2}^2+\ldots+a_1^2\cdot\ldots\cdot a_{[k/2]}^2}.\nonumber
\end{eqnarray}
On the other hand
\begin{equation}\nonumber
\frac{1}{l!}j^*\alpha^l=\sum_{j_1<\ldots<j_l}a_{j_1}\cdot\ldots\cdot a_{j_l}e_{2j_1-1}\wedge e_{2j_1}\wedge\ldots\wedge e_{2j_l-1}\wedge e_{2j_l},
\end{equation}
so that
\begin{equation}\nonumber
g(\frac{1}{l!}j^*\alpha\wedge\star \varpi_L,\frac{1}{l!}j^*\alpha\wedge\star \varpi_L)=\sum_{j_1<\ldots<j_l}a_{j_1}^2\ldots a_{j_l}^2,\quad 2l\leq k.
\end{equation}
Summing and taking the square root yields precisely~(\ref{determinant}).  
\end{prf}

Next we show that in analogy with the classical counterpart discussed in Section~\ref{classcal}, the spinors inducing a $G_L\times G_R$--structure define calibrations in the generalised sense. In particular, we will recover ``classical'' calibrated $k$--planes as a special case. First, we associate with any isotropic pair $(L,F)$ an isometry $\mc{R}_{L,F}$, the {\em gluing matrix} (this jargon stems from physics, cf. Section~\ref{wspov}). To start with, let $B=0$ and consider the isotropic subspace $W_L$. We can think of this space as the graph of the anti--isometry $P_L:D^+\to D^-$ induced by the orthogonal splitting $T\oplus T^*=D^+\oplus D^-$ with isometries $\pi_{0\pm}:X\in(T,\pm g)\mapsto X^{\pm}\in(D^{\pm},g_{\pm})$ (cf. Section~\ref{genmetrics}). Indeed, we have
\begin{equation}\nonumber
0=(X\oplus P_LX,X\oplus P_LX)=(X,X)+(P_LX,P_LX)=g_+(X,X)+g_-(P_LX,P_LX).
\end{equation}
If we choose an adapted orthonormal basis $e_1,\ldots,e_k\in L$, $e_{k+1},\ldots e_n\in L^{\perp}$, then the matrix of $P_L$ associated with the basis $d^{\pm}_l=\pi_{0\pm}(e_l)=e_l\oplus\pm e^l$ of $D^{\pm}$ is
\begin{equation}\label{pumatrix}
P_L=\left(\begin{array}{cc} Id_k & 0\\ 0 & -Id_{n-k}\end{array}\right).
\end{equation}
Pulling this operator back to $T$ via $\pi_{0\pm}$ yields the isometry
\begin{equation}\label{rumatrix}
\mc{R}_L=\pi_{0-}^{-1}\circ P_L\circ\pi_{0+}:(T,g)\to (T,g).
\end{equation}
Its matrix representation with respect to an orthonormal basis adapted to $L$ is just~(\ref{pumatrix}). In the general case, $W_{L,F}$ is the graph of the anti--isometry $P_{L,F}:V^+\to V^-$ with associated gluing matrix
\begin{equation}\nonumber
\mc{R}_{L,F}
=\pi^{-1}_-\circ P_{L,F}\circ\pi_+,
\end{equation}
where now $\pi_{\pm}:(T,\pm g)\to V^{\pm}=e^{2B}D^{\pm}$. Let $\mc{F}=F-j^*B$. Changing, if necessary, the orthonormal basis on $L$ such that $\mc{F}=\sum_{l=1}^{[k/2]}f_le_{2l-1}\wedge e_{2j}$, 
\begin{equation}\nonumber
e_1\oplus f_1e^2,e_2\oplus-f_1e^1,\ldots,e^{k+1},\ldots,e^n
\end{equation}
is a basis of $e^{2B}W_{L,F}$ by~(\ref{isorep}). Decomposing the first $k$ basis vectors into the $D^{\pm}$--basis $d^{\pm}_j$ yields
\begin{equation}\nonumber
\begin{array}{rcr}
2(e_{2l-1}\oplus f_l e^{2l}) & = & d_{2l-1}^++f_l d_{2l}^+\oplus d_{2l-1}^--f_l d_{2l}^-\\
2(e_{2l}\oplus -f_l e^{2l-1}) & = & -f_l d_{2l-1}^++d_{2l}^+\oplus f_l d_{2l-1}^-+d_{2l}^-,
\end{array}\,,\quad l\leq k
\end{equation}
while $2e^l=d_l^+\oplus-d_l^-$ for $l=k+1,\ldots,n$.
Written in the new $D^{\pm}$--basis 
\begin{equation}\nonumber
\begin{array}{ll}
w_{2l-1}^+=d_{2l-1}^++f_ld_{2l}^+,\, w_{2l}^+=-f_ld_{2l-1}^++d_{2l}^+,\,l\leq k,\quad & w_l^+=d_l^+,\,l>p\mbox{ on }D^+\\ w_{2l-1}^-=d_{2l-1}^--f_ld_{2l}^-,\, w_{2l}^-=\phantom{-}f_ld_{2l-1}^-+d_{2l}^-,\,l\leq k,\quad & w_l^-=d_l^-,\, l>k\mbox{ on }D^-,
\end{array}
\end{equation} 
the matrix $\mc{R}_{L,F}$ is just~(\ref{pumatrix}).
The base of change matrix for $d_l^+\to w_l^+$ is given by the block matrix $A=(A_1,\ldots,A_k,id_{n-k})$, where for $(d_{2l-1}^+,d_{2l}^+)\to (w_{2l-1}^+,w_{2l}^+)$, $l\leq k$,
\begin{equation}\nonumber
A_l=\left(\begin{array}{rr} 1 & f_l\\ -f_l & 1\end{array}\right).
\end{equation}
The basis change $w_l^-\to d_l^-$ is implemented by $B=(B_1,\ldots,B_k,id_{n-k})$, where for $(d_{2l-1}^-,d_{2l}^-)\to (w_{2l-1}^-,w_{2l}^-)$, $l\leq k$,
\begin{equation}\nonumber
B_l=\left(\begin{array}{rr} 1 & -f_l\\ f_l & 1\end{array}\right).
\end{equation}
Computing $B\circ(Id_k,-Id_{n-k})\circ A^{-1}$ and pulling back to $T$ via $\pi_{\pm}$, we finally find
\begin{equation}\label{ruf}
\mc{R}_{L,F}=\left(\begin{array}{cc} \big(j^*(g-B)+F\big)\big(j^*(g+B)-F\big)^{-1} & 0\\ 0 & -Id_{n-k}\end{array}\right)
\end{equation}
with respect to some orthonormal basis adapted to $L$.

\begin{prp}\label{liesinspin}
For any isotropic pair $(L,F)$, the element 
\begin{equation}\nonumber
\mf{J}^{-1}\big(e^{-B}\bullet\frac{\tau_{L,F}}{\eunorm{\tau_{L,F}}_{\mc{G}}}\big)\in\cliff(T,g)
\end{equation} 
lies in $Pin(T,g)$. Moreover, its projection to $O(T,g)$ equals the gluing matrix $\mc{R}_{L,\mc{F}}$. 
\end{prp}

\begin{prf}
Again let $e_1,\ldots,e_n$ be an adapted orthonormal basis so that $\mc{F}=\sum_{l=1}^{[k/2]}f_le_{2l-1}\wedge e_{2l}$. Applying a trick from~\cite{bkop97}, we define
\begin{eqnarray*}
\arctan{\widetilde{\mc{F}}} & = & \sum_l\arctan(f_l)e_{2l-1}\cdot e_{2l}\\
& = & \frac{1}{2i}\sum_l\ln\frac{1+if_l}{1-if_l}e_{2l-1}\cdot e_{2l}\in\mf{spin}(n)\subset\cliff(T,g)
\end{eqnarray*}
and show that $\mf{J}\big(\exp(\arctan\widetilde{\mc{F}})\cdot\widehat{\star \varpi_L}\big)=e^{-B}\bullet\tau_{L,F}/\!\!\eunorm{\tau_{L,F}}_{\mc{G}}\,\in\Lambda^*T^*$, where $\exp$ takes values in $Spin(T,g)$. Since the elements $e_{2l-1}\cdot e_{2l}$, $e_{2m-1}\cdot e_{2m}$ commute, exponentiation yields
\begin{eqnarray*}
e^{\arctan{\widetilde{\mc{F}}}} & = & \prod_m e^{\arctan(f_k)e_{2m-1}\cdot e_{2m}}\\
& = & \prod_m\big(\cos\big(\arctan(f_m)\big)+\sin\big(\arctan(f_m)\big)e_{2m-1}\cdot e_{2m}\big)\\
& = & \prod_m\big(\frac{1}{\sqrt{1+f_m^2}}+\frac{f_m}{\sqrt{1+f_m^2}}e_{2m-1}\cdot e_{2m}\big)\\
& = & \frac{\big(1+\sum_lf_le_{2l-1}\cdot e_{2l}+\sum_{l<r}f_l\cdot f_re_{2l-1}\cdot e_{2l}\cdot e_{2r-1}\cdot e_{2r}+\ldots\big)}{\prod_m(\sqrt{1+f_m^2})}\\
& = &
\frac{1}{\sqrt{\det(j^*g-\mc{F})}}\mf{J}^{-1}(1+\mc{F}+\frac{1}{2}\mc{F}\wedge\mc{F}+\ldots),
\end{eqnarray*}
where we used the classical identities $\cos\arctan x=1/\sqrt{1+x^2}$, $\sin\arctan x=x/\sqrt{1+x^2}$. By Lemma~\ref{det}, we finally get
\begin{equation}\nonumber
\mf{J}(e^{\arctan{\widetilde{\mc{F}}}}\cdot\widehat{\star \varpi_L})=\mf{J}(e^{\arctan\widetilde{\mc{F}}})\wedge \mf{J}(\widehat{\star \varpi_L})=e^{-B}\bullet \frac{\tau_{L,F}}{\eunorm{\tau_{L,F}}_{\mc{G}}}.
\end{equation}
The projection $\pi_0:Pin(T,g)\to O(T,g)$ via $\pi_0$ gives indeed the induced gluing matrix: Firstly,
\begin{equation}\nonumber
\pi_0\big(e^{\arctan\widetilde{F}}\cdot \widehat{\star\varpi_L}\big)=e_{SO(T)}^{\pi_{0*}(\arctan\widetilde{F})}\circ\pi_0(\widehat{\star\varpi_L}).
\end{equation}
Now 
\begin{eqnarray*}
e^{\pi_{0*}\big(\arctan(f_l)e_{2l-1}\cdot e_{2l}\big)} & = &
e^{2\arctan(f_l)e_{2l-1}\wedge e_{2l}}\\
& = & \cos\big(2\arctan(f_l)\big)+\sin\big(2\arctan(f_l)\big)e_{2l-1}\wedge e_{2l}\\
& = & \frac{1-f_l^2}{1+f_l^2}+\frac{2f_l}{1+f_l^2}e_{2l-1}\wedge e_{2l}\\
& = & \frac{1}{1+f^2_l}\big(\left(\begin{array}{cc} 1-f^2_l & 0\\ 0 & 1-f^2_l\end{array}\right)+\left(\begin{array}{cc} 0 & -2f_l\\ 2f_l & 0\end{array}\right)\big)
\end{eqnarray*}
which yields the matrix $(j^*g+\mc{F})(j^*g-\mc{F})^{-1}$, while the block $-Id_{n-k}$ in the gluing matrix is accounted for by the projection of the volume form.
\end{prf}

\noindent The fact that $e^{-B}\bullet\tau_{L,F}/\!\eunorm{\tau_{L,F}}_{\mc{G}}$ can be identified with an element of $Pin(T,g)$ enables us to prove the following

\begin{thm}\label{spincrit}
Let $\Psi_L,\Psi_R$ be two chiral unit spinors of the complex $Spin(T)$--representation $\Delta_n$. The $T\oplus T^*$--spinors $\rho^{ev,od}=e^B\bullet\Re\big([\Psi_L\otimes\Psi_R]^{ev,od}\big)$ satisfy $\langle\rho,\tau_{L,F}\rangle\leq\eunorm{\tau_{L,F}}_{\mc{G}}$. Moreover, an isotropic pair $(L,F)$ meets the boundary if and only if
\begin{equation}\nonumber
\mc{A}(\Psi_L)=\pm(-1)^{m(m+1)/2+k}i^m(e^{-B}\bullet\frac{\tau_{L,F}}{\eunorm{\tau_{L,F}}_{\mc{G}}})\cdot\Psi_R.
\end{equation} 
for $n=2m$ and $\Psi_R\in\Delta_{\pm}$ and
\begin{equation}\nonumber
\mc{A}(\Psi_L)=(-1)^{m(m+1)/2}i^{m+1}(e^{-B}\bullet \frac{\tau_{L,F}}{\eunorm{\tau_{L,F}}_{\mc{G}}})\cdot\Psi_R.
\end{equation} 
for $n=2m+1$. In particular, $\rho^{ev,od}$ define a calibration.
\end{thm}

\begin{prf} 
Since $e^{\mc{F}}\wedge\widehat{\star \varpi_L}=(\star e^{\mc{F}}\llcorner \varpi_L)^{\wedge}$, we obtain
\begin{eqnarray}
\langle e^B\bullet\Re [\Psi_L\otimes\Psi_R]^{ev,od},\tau_{L,F}\rangle
& = & \langle \Re[\Psi_L\otimes\Psi_R]^{ev,od},e^{\mc{F}}\wedge\widehat{\star\varpi_L}\rangle\nonumber\\
& = &  g(\Re[\Psi_L\otimes\Psi_R]^{ev,od},e^{\mc{F}}\llcorner\varpi_L)\nonumber\\
& = & 
\sum\Re q(\mc{A}(\Psi_L),e_I\cdot\Psi_R)g(e_I,e^{\mc{F}}\llcorner\varpi_L)\nonumber\\ 
& = & 
\Re q(\mc{A}(\Psi_L),e^{\mc{F}}\llcorner\varpi_L\cdot\Psi_R).\label{cauchy} 
\end{eqnarray} 
On the other hand,
\begin{eqnarray*}
e^{\mc{F}}\llcorner\varpi_L\cdot\Psi_R & = & (-1)^{k(n-k)}(e^{\mc{F}}\llcorner\star\star\varpi_L)\cdot\Psi_R\\
& = & (-1)^{k(n-k)}\star(e^{-\mc{F}}\wedge\star\varpi_L)\cdot\Psi_R\\
& = & (-1)^{k(n-k)}(e^{\mc{F}}\wedge\widehat{\star\varpi_L})\cdot\varpi_g\cdot\Psi_R\\
& = & (-1)^{k(n-k)}\sqrt{\det(j^*g-\mc{F})}(e^{-B}\bullet\tau_{L,F})\cdot\varpi_g \cdot\Psi_R.
\end{eqnarray*}
Recall from~(\ref{volaction}) the action of $\varpi_g$ on $\Delta_{\pm}$. Since $e^{-B}\bullet\tau_{L,F}/\!\eunorm{\tau_{L,F}}_{\mc{G}}\in Spin(T,g)$, we get
\begin{equation}\nonumber
\eunorm{e^{\mc{F}}\llcorner\varpi_L \cdot\Psi_R}_q=\sqrt{\det(j^*g-\mc{F})}\eunorm{e^{-B}\bullet\frac{\tau_{L,F}}{\eunorm{\tau_{L,F}}_{\mc{G}}}}_g\cdot\eunorm{\Psi_R}_q=\sqrt{\det(j^*g-\mc{F})}
\end{equation}
(where the norm is taken for the Hermitian inner product $q$ on $\Delta_n$ and $g$ on $T$ respectively). As a consequence,~(\ref{cauchy}) is less than or equal to $\eunorm{\tau_{L,F}}_{\mc{G}}$ by the Cauchy--Schwarz inequality for the norm $\eunorm{\cdot}_q$. Moreover, equality holds precisely if $\mc{A}(\Psi_L)=(-1)^{k(n-k)}(e^{-B}\bullet\tau_{L,F})\cdot\varpi_g\cdot\Psi_R$. As there always exists a subspace $L$ such that $\mc{A}(\Psi_R)=\varpi_L\cdot\Psi_L$, we can choose $(L,j^*B)$ as a calibrated isotropic pair. Hence, the spinor $\rho^{ev,od}=e^B\bullet\Re\big([\Psi_L\otimes\Psi_R]^{ev,od}\big)$ defines a calibration. 
\end{prf}

\begin{ex}
Consider a $G_2$-- or $Spin(7)$--structure on $T$ induced by the real unit spinor $\Psi$. Then $(L,0)$ defines a calibrated plane for $\Re\big([\Psi\otimes\Psi]^{ev,od}\big)$ if and only if
\begin{equation}\nonumber
\Psi=\star\varpi_L\cdot\Psi=\varpi_L\cdot\Psi
\end{equation}
as $\varpi_g$ acts as the identity on $\Delta_n$ for $n=7,\,8$. This is precisely the condition found by Dadok and Harvey (cf. Section~\ref{classcal}). In general, if $L$ is a calibrated plane with respect to some classical calibration form arising as the degree $k$ component of $\rho=\Re\big([\Psi\otimes\Phi]^{ev,od}\big)$, the isotropic pair $(L,0)$ is calibrated with respect to $\rho$ in the sense of Definition~\ref{caldef} with $B$--field transform $(L,j^*B)$. More examples will be given in Sections~\ref{calsub} and~\ref{buscherrules}. In particular, we will see that even in a ``straight'' setting with $B=0$, Definition~\ref{caldef} is more general as we can have calibrated isotropic pairs $(L,F)$ with non--trivial $F$.
\end{ex}
%
\subsection{Generalised calibrated submanifolds}
\label{calsub}
We introduce the notion of a generalised calibrated submanifold next. Let $M$ be a smooth manifold with a generalised Riemannian metric $V^+\subset\mb{E}$ corresponding to $(g,H)$.

\begin{definition}\hfill\newline
\noindent{\rm (i)} A spinor field $\rho\in\Gamma\big(\mb{S}(\mb{E})_{\pm}\big)$ is called a {\em generalised calibration} if and only if $\rho_p\in\mb{S}(\mb{E})_{\pm\,p}$ defines a generalised calibration for every $p\in M$.

\noindent{\rm (ii)} An {\em isotropic pair} $(L,\mc{F})$ for $(M,g,H)$ consists of an embedded oriented submanifold $j=j_L:L\hookrightarrow M$, together with a $2$--form $\mc{F}\in\Omega^2(L)$ such that
\begin{equation}\label{gauge}
d\mc{F}+j^*H=0.
\end{equation}

\noindent{\rm (iii)} An isotropic pair $(L,\mc{F})$ is said to be {\em calibrated by} $\rho$ if and only if for every point $p\in j(L)\subset M$ and for one (and thus for any) $\tau_p$ in the half--line of $\mb{S}(\mb{E})_{\pm\,p}$ induced by the pure spinor field
\begin{equation}\nonumber
\tau_{L,\mc{F}}=\mf{L}^{-1}(e^{\mc{F}}\wedge\widehat{\star\varpi_L})\in\Gamma\big(\mb{S}(\mb{E})_{\pm|L}\big),
\end{equation}
we have $\langle\rho_p,\tau_{L,\mc{F}}\rangle=\eunorm{\tau_p}_{\gtilde}$. Here, we think of $\mc{F}$ as a section of $\Lambda^2T^*\!M_{|j(L)}$ via the pull--back induced by orthogonal projection $p_L:TM_{|j(L)}\to Tj(L)$.
\end{definition}

\noindent Again, we will drop the qualifier ``generalised'' and simply speak of a calibration form.  Also, we will usually think of $L$ as a subset of $M$ and identify $L$ with $j(L)$.

An easy corollary is the following ``form criterion'' for an isotropic pair to be calibrated.

\begin{prp}\label{fromcritmanif}
An isotropic pair $(L,\mc{F})$ for $(M,g,H)$ with $k=\dim L$ is calibrated if and only if
\begin{equation}\nonumber
g(e^{-\mc{F}}\wedge \mf{L}(\rho),\varpi_L)\leq\eunorm{\tau_{L,\mc{F}}}_{\mc{G}}
\end{equation}
or equivalently,
\begin{equation}\nonumber
[e^{-\mc{F}}\wedge j^*_L\mf{L}(\rho)]^k\leq\sqrt{\det(j^*_Lg-\mc{F})}\varpi_L.
\end{equation}
\end{prp}

\noindent This follows from
\begin{equation}\nonumber
\langle\rho,\tau_{L,\mc{F}}\rangle=\star[\mf{L}(\rho)\wedge e^{-\mc{F}}\wedge\star \varpi_L]^n=g(e^{-\mc{F}}\wedge \mf{L}(\rho),\varpi_L),
\end{equation}
while
\begin{equation}\nonumber
\eunorm{\tau_{L,\mc{F}}}^2_{\mc{G}}=\star[\mf{L}(\tau_{L,\mc{F}})\wedge\widehat{\gtilde\big(\mf{L}(\tau_{L,\mc{F}})\big)}]^n=g(e^{\mc{F}}\wedge\widehat{\star \varpi_L},e^{\mc{F}}\wedge\widehat{\star \varpi_L}).
\end{equation}

Calibrated isotropic pairs minimise an extension of the volume functional, namely the {\em brane energy}
\begin{equation}\nonumber
\mc{E}_{\phi,\gamma}(L,\mc{F})=\int_Le^{-\phi}\eunorm{\tau_{L,\mc{F}}}_{\mc{G}}-\int_Le^{-\phi}j^*_L[e^{-\mc{F}}\wedge \mf{L}(\gamma)],
\end{equation}
where $\gamma$ is a spinor of suitable parity and $\phi\in C^{\infty}(M)$ a dilaton field. 
The terminology comes from the local description of the spinor norm $\eunorm{\tau_{L,\mc{F}}}_{\mc{G}}$. Let $C=\mf{L}(\gamma)$. These are the so--called {\em Ramond--Ramond potentials} (cf. Section~\ref{physics2}). Over $U_a$, $H_{|U_a}=dB^{(a)}_a$, so that for a disk $D\subset U_a\cap L$, the brane energy is just the sum of the two terms
\begin{equation}\nonumber
\int_D e^{-\phi}\eunorm{\tau_{L,\mc{F}}}_{\mc{G}}\,=\int_D\sqrt{\det\big(j^*_L(g-B^{(a)}_a)+F^{(a)}\big)}\varpi_D
\end{equation}
and 
\begin{equation}\nonumber
\int_De^{-\phi}j^*_L[e^{-\mc{F}}\wedge \mf{L}(\gamma)]=\int_De^{-\phi}j^*_L[e^{-F^{(a)}+B^{(a)}}\wedge C]
\end{equation}
In physics, these integrals are known as the {\em Dirac--Born--Infeld} and the {\em Wess--Zumino} term of the D--brane energy (cf.~(\ref{phys_dbi}) and~(\ref{phys_wz})).

\begin{definition}\hfill\newline
\noindent{\rm (i)} Two isotropic pairs $(L,\mc{F})$ and $(L',\mc{F}')$ of dimension $k$ are said to be {\em homologous}, if there exists an isotropic pair $(\widehat{L},\widehat{\mc{F}})$ of dimension $k+1$ with $\partial\widehat{L}=L-L'$, $\widehat{\mc{F}}_{|L}=\mc{F}$ and $\widehat{\mc{F}}_{|L'}=\mc{F}'$. 

\noindent{\rm (ii)} We call an isotropic pair $(L,\mc{F})$ {\em locally brane--energy minimising for} $\mc{E}_{\phi,\gamma}$ if
\begin{equation}\nonumber
\mc{E}_{\phi,\gamma}(D,\mc{F}_{|D})\leq\mc{E}_{\phi,\gamma}(D',\mc{F}')
\end{equation}
for any embedded disc $D\subset L$ and homologous isotropic pair $(D',\mc{F}')$ (with $D'$ not necessarily embedded into $L$).
\end{definition}

\begin{thm}\label{minimum}
Let $\rho$ be a calibration with
\begin{equation}\label{intcond}
d_{\phi}\rho=d_{\phi}\gamma,\mbox{ i.e.}\quad d_He^{-\phi}\mf{L}(\rho)=d_He^{-\phi}\mf{L}(\gamma).
\end{equation}
Then any calibrated isotropic pair $(L,\mc{F})$ is locally brane--energy minimising for $\mc{E}_{\phi,\gamma}$.
\end{thm}

\begin{prf}
With the notation of the previous definition, let $(\widehat{D},\widehat{\mc{F}})$ be an isotropic pair such that $\partial\widehat{D}=D-D'$, $\widehat{\mc{F}}_{|D}=\mc{F}$ and $\widehat{\mc{F}}_{|D'}=\mc{F}'$. By Stokes' theorem
\begin{eqnarray*}
0 & = & \int_{\widehat{D}}j^*_{\widehat{D}}\big[e^{-\widehat{\mc{F}}}\wedge d_He^{-\phi}\mf{L}(\rho-\gamma)\big]\\
& = & \int_{\widehat{D}}j^*_{\widehat{D}}\big[d\big(e^{-\widehat{\mc{F}}}\wedge e^{-\phi}\mf{L}(\rho-\gamma)\big)\big]\\
& = & \int_{D-D'}j^*_{\partial\widehat{D}}[e^{-\widehat{\mc{F}}_{|\partial\widehat{D}}}\wedge e^{-\phi}\mf{L}(\rho-\gamma)],
\end{eqnarray*}
so that
\begin{equation}\nonumber
\int_De^{-\phi}j^*_D\big[e^{-\mc{F}}\wedge\mf{L}(\rho-\gamma)\big]=\int_{D'}e^{-\phi}j^*_{D'}\big[e^{-\mc{F}'}\wedge\mf{L}(\rho-\gamma)\big].
\end{equation}
Since $D\subset L$ and $(L,\mc{F})$ is calibrated,
\begin{equation}\nonumber
\int_Dj^*_De^{-\phi}[e^{-\mc{F}}\wedge\mf{L}(\rho)]=\int_De^{-\phi}\eunorm{\tau_{L,\mc{F}}}_{\mc{G}}
\end{equation} 
by the previous lemma, while for $D'$, the integral over the norm is greater or equal than $\int_{D'}j^*_{D'}[e^{-\phi}e^{-\mc{F}'}\wedge \mf{L}(\rho)]$.
\end{prf}

\begin{rmk}
In~\cite{jewi05b}, the equation $d_{\phi}\rho=d_{\phi}\gamma$ was shown to describe type II string compactification on $6$ or $7$ dimensions which are governed by the {\em democratic formulation} of Bergshoeff et al~\cite{bkop01}.
\end{rmk}

Examples of calibrations are provided by $G_L\times G_R$--structures: The following proposition is an immediate consequence of Theorem~\ref{spincrit}.

\begin{prp}
If $\rho\in\Gamma\big(\mb{S(E)}_{\pm}\big)$ is a spinor field such that $\mf{L}(\rho)=[\Psi_L\otimes\Psi_R]$ for (chiral) unit spinor fields $\Psi_{L,R}\in\Delta(TM)$, then $\Re(\rho)$ defines a calibration form. Moreover, an isotropic pair is calibrated if and only if 
\begin{equation}\label{susybreak1}
\mc{A}(\Psi_L)=\pm(-1)^{m(m+1)/2+k}i^m\frac{\mf{L}(\tau_{L,\mc{F}})}{\eunorm{\mf{L}(\tau_{L,\mc{F}})}_g}\cdot\Psi_R.
\end{equation} 
for $n=2m$ and $\Psi_R\in\Delta_{\pm}$, and
\begin{equation}\label{susybreak2}
\mc{A}(\Psi_L)=(-1)^{m(m+1)/2}i^{m+1}\frac{\mf{L}(\tau_{L,\mc{F}})}{\eunorm{\mf{L}(\tau_{L,\mc{F}})}_g}\cdot\Psi_R.
\end{equation} 
for $n=2m+1$, where $\eunorm{\mf{L}(\tau_{L,\mc{F}})}_g=\sqrt{g\big(\mf{L}(\tau_{L,\mc{F}}),\mf{L}(\tau_{L,\mc{F}})\big)}$.
\end{prp}

\begin{rmk}
Equations~(\ref{susybreak1}) and~(\ref{susybreak2}) reflect the physical fact that calibrated pairs ``break the supersymmetry'': The spinor field $\tau_{L,\mc{F}}$ induces a linear relation between $\Psi_L$ and $\Psi_R$ (or the {\em supersymmetry parameters} in physicists's language, cf. Section~\ref{wspov}), so that over $L$ we are left with only one independent supersymmetry parameter. On the tangent space level, the two $G$--structures $P_{G_{L,R}}$ inside the orthonormal frame bundle are related by the corresponding gluing matrix.
\end{rmk}

\paragraph{Examples.}
By the previous proposition, the $\mb{E}$--spinor fields which induce a generalised $SU(3)$--, $G_2$-- or $Spin(7)$--structure, define calibrations. For straight generalised $SU(3)$--, $G_2$-- and $Spin(7)$--structures (cf. Section~\ref{gengstruc}), we recover some well--known examples coming from physics which involve a non--trivial $\mc{F}$. We will construct further examples in Section~\ref{tduality} using the device of T--duality.

\bigskip

\noindent(i) {\em Straight $SU(3)$--structures:} Let $(M^6,g,\Psi)$ be a classical $SU(3)$--structure whose induced straight $SU(3)\times SU(3)$--structure is given by $\mf{L}(\rho_0)=e^{-i\omega}$ and $\mf{L}(\rho_1)=\Omega$. Starting with the even spinor field, we have $\Re\big(\mf{L}(\rho^{ev})\big)=1-\omega^2/2$. By Proposition~\ref{fromcritmanif}, we find as calibration condition
\begin{equation}\label{su3_even}
\frac{1}{k!}\left(j_L^*\omega+\mc{F}\right)^k = \sqrt{\det(j_L^*g-\mc{F})}\varpi_L,
\end{equation}
which defines so--called {\em B--branes}. These wrap holomorphic cycles in accordance with results from~\cite{kali04},~\cite{mmms00}. Calibrating with respect to $\Re\big(\mf{L}(\rho^{od})\big)=\Re(\Omega)$ yields so--called {\em A--branes} defined by
\begin{equation}\nonumber
\left[e^{-\mc{F}}\wedge j_L^*\Re(\Omega)\right]^k = \sqrt{\det(j_L^*g-\mc{F})}\varpi_L.
\end{equation}
The calibrated pairs are therefore odd--dimensional. For $k=3$,  $j_L^*\Re\Omega=\varpi_L$ is the condition for a special Lagrangian cycle. For $k=5$ we need a non--vanishing $\mc{F}$ and obtain
\begin{equation}\nonumber
j_L^*\Re\Omega\wedge \mc{F}=\sqrt{\det(j_L^*g-\mc{F})}\varpi_L
\end{equation}
which is the condition found in~\cite{kali04} for a {\em coisotropic} A--brane. In the non--straight case, $\Re\mf{L}(\rho^{od})$ can also contain a 1-- and 5--form part. For a related discussion on this aspect see Section 4 of~\cite{begr06}.

\medskip

\noindent(ii) {\em Straight $G_2\!$--structures:}
Let $(M^7,g,\Psi)$ be a classical $G_2\!$--structure whose induced straight $G_2\times G_2$--structure is given by
$\mf{L}(\rho^{ev})=1-\star\varphi$. The calibration condition reads
\begin{equation}\nonumber
e^{-\mc{F}}\wedge j_L^*(1-\star\varphi)\leq\sqrt{\det(j_L^*g-\mc{F})}\varpi_L.
\end{equation}
A coassociative submanifold satisfies $j_L^*\star\!\varphi=\varpi_L$ and is therefore calibrated if $\mc{F}=0$. For non--trivial $\mc{F}$, we find $\mc{F}\wedge \mc{F}/2-j^*_L\star\varphi=\sqrt{\det(j_L^*g-\mc{F})}\varpi_L$. Now $\mc{F}\wedge \mc{F}/2=\mathrm{Pf}(\mc{F})\varpi_L$ and $\det(j_L^*g-\mc{F})=1-\mathrm{Tr}(\mc{F}^2)/2+\det(\mc{F})$. Squaring this shows the equality to hold for an isotropic pair $(L,\mc{F})$, where $L$ is a coassociative submanifold and $\mc{F}$ a closed $2$--form with $2\mathrm{Pf(\mc{F})}=-\mathrm{Tr}(\mc{F}^2)/2$, that is, $\mc{F}$ is anti--self--dual (cf.~\cite{mmms00}).

\medskip

\noindent(iii) {\em Straight $Spin(7)$--structures:}
Let $(M^8,g,\Psi)$ be a classical $Spin(7)$--structure whose induced straight $Spin(7)\times Spin(7)$--structure is given by $\mf{L}(\rho)=1-\Omega+\varpi_g$. Again, Cayley submanifolds are calibrated for $\mc{F}=0$. Further, as in (ii), they remain so if turning on an anti--self--dual gauge field $\mc{F}$ (cf.~\cite{mmms00}).
%
%
%
%
%
\section{Branes in string theory}
\label{physics}
We briefly interlude to outline how in type II string theory the calibration condition as presented above arises from both the world--sheet and the target space point of view. For a detailed introduction see for instance~\cite{jo03} or~\cite{po96}. The material of this section is independent of the mainstream development of the paper.

Generally speaking, D--branes arise as boundary states in the conformal field theory on the string worldsheet. In the supergravity limit, ignoring corrections of higher order in $\alpha'$ (the {\em string tension}), they can be described as submanifolds in the ten dimensional target space which extremalise a certain energy functional. In the special case of type II string compactifications on Calabi--Yau manifolds, D--branes can be classified by the derived category of coherent sheaves and the Fukaya category for type IIB and IIA respectively. Both notions are connected by Mirror symmetry, which in the case of toroidal fibrations ought to be realised through T--duality~\cite{syz96}.
%
\subsection{World--sheet point of view}
\label{wspov}
We start with a two--dimensional conformal field theory on the string
world--sheet $\Sigma$, parametrised by a space--like coordinate $s$ and
a time--like coordinate $t$.
D--branes arise when considering open string solutions (i.e. where the string
is homeomorphic to an open interval) with Dirichlet boundary conditions.

Let us first consider the simplest case without background fields, delaying this issue to the end of this section. Varying the world--sheet action on $\Sigma$
(which depends on the embedding functions $X^\mu$, $\mu=0,\ldots,9$, taken as
coordinates of the $10$--dimensional target space), we find the boundary term
\begin{equation}\nonumber
\delta I_{\mathrm{boundary}}=\int_{\partial\Sigma}\partial_s X_\mu\delta X^\mu.
\end{equation}
There are two kinds of solutions at the boundary of the worldsheet ($s=\{0,\pi\}$), namely
\begin{equation}\label{eq_bc}
\begin{array}{lcl}
\partial_s X^\mu|_{s=0,\pi} = 0&\quad&\mbox{(von~Neumann boundary condition)}\\\vspace{-10pt} & & \\&\mbox{or}&\\
\delta X^\mu = 0 \Rightarrow X^\mu|_{s=0,\pi}=\mbox{const.}&\quad&\mbox{(Dirichlet boundary condition).}
\end{array}
\end{equation}
Choosing the von~Neumann boundary condition for $\mu=0,\ldots,p$ and the Dirichlet boundary condition for
$\mu=p+1,\ldots,d$, we define open strings whose endpoints can move along $p$ spacial
dimensions (the zeroth coordinate parametrising time), thus sweeping out a $p+1$--dimensional surface, a {\em D(p)--brane}.

Since we are dealing with a world--sheet theory that is supersymmetric,
the boundary conditions for the world--sheet fermions need to be taken into account, too. Without D--branes we start with $(1,1)$ worldsheet supersymmetry,
generated by two spinorial parameters, $\epsilon_{L,R}$. The supersymmetric partners of the coordinate functions $X^\mu$ are represented by the set of worldsheet fermions $\psi^\mu_{L,R}$. The boundary conditions for these read
\begin{equation}\nonumber
\big(-\psi_{L\mu}\delta\psi_L^\mu+\psi_{R\mu}\delta\psi_R^\mu\big)|_{s=0,\pi}=0,
\end{equation}
which can be solved by taking $\psi_{L\mu}=\pm\psi_{R\mu}$. However, the fermionic boundary conditions are related to their bosonic counterpart by supersymmetry, so only one choice is effectively possible. Moreover, the supersymmetry transformations between fermions and bosons imply that the supersymmetry parameters are related by $\epsilon_R=-\epsilon_L$. Hence we are left with only one linearly independent supersymmetry parameter so that worldsheet supersymmetry is broken to $(1,0)$. By convention, a $+$ occurs if the coordinate $\mu$ satisfies the von~Neumann boundary condition. For directions satisfying the Dirichlet boundary conditions we then obtain a minus sign, i.e. $\psi_{L\mu}=-\psi_{R\mu}$.
These boundary conditions can be suitably encoded in the {\em gluing matrix} $\mc{R}$
(\ref{rumatrix}), relating
$\psi_L$ and $\psi_R$ via $\psi_L=\mc{R}\psi_R$, $\mc{R}$ being the diagonal matrix (with respect to the coordinates $\mu=0,\ldots 9$) defined by $+1$
in von~Neumann and $-1$ in Dirichlet directions.

As a consequence of the D--brane boundary conditions, the symmetry group on the branes
reduces to $SO(1,p)$. After quantisation, the zero mode spectrum of the open strings contains a vector multiplet of $\mc{N}=1$ supersymmetry, that is, an element of the $SO(p-1)$ vector representation plus its fermionic counterpart. The vector multiplet also transforms under a $U(1)$ gauge symmetry. In particular, we obtain a gauge field $F$ -- the curvature 2--form coming from the $U(1)$--connection.
More generally, we can take $N$ branes on top of each other to obtain a $U(N)$
gauge group, reflecting the fact that we have to include the fields from strings
stretching between different branes of one stack. In any case, including the field strength leads to a generalised gluing condition of the form (\ref{ruf}).
To see how this happens, consider again the boundary conditions~(\ref{eq_bc}). First note that the Dirichlet condition is equivalent to
$\partial_t X^\mu=0$. Including a nontrivial gauge field, or more generally a two--form field--strength background $\mc{F}$, which also includes the background $B$--field (cf. section~\ref{physics2}), changes the boundary conditions into mixed Dirichlet-- and von~Neumann--type as follows,
\begin{equation}\nonumber
\partial_s X|_{\partial L}+g^{-1}(\partial_t X\llcorner\mc{F})|_{\partial L}=0,
\end{equation}
where $g$ is the space--time metric and restriction to $\partial L$ indicates restriction to the boundary of the D--brane.
%
\subsection{Target space point of view}
\label{physics2}
For type II supergravity
theory on $M$, that is, at the zero--mode level of the corresponding string theory, we have the following field content in the so--called {\em closed string sector}:
\begin{itemize}
\item a metric $g$, a 2--form $B$ (the {\em B--field}) and a scalar function $\phi$ (the {\em dilaton}) in the NS-NS (NS=Neveu--Schwarz) sector.
\item a mixed form $C$ (the {\em R-R--potentials}) in the R-R (R=Ramond) sector, which is odd or even in type~IIA or~IIB respectively. 
\end{itemize}
In the limit of supergravity, D--branes become non--perturbative objects and belong thus to the background geometry where they can be conceived as submanifolds of the $10$--dimensional target manifold $M$. In this picture, $H$ is given locally as the field strength $H=dB$.

On the other hand, the open string sector leads to an effective field theory on the world--volume of a D--brane. If $L$ denotes the corresponding submanifold, this field theory is given by the {\em Dirac--Born--Infeld action}
\begin{equation}\label{phys_dbi}
I_\mathrm{DBI}(L,F)=-N \int_L e^{-\phi}\sqrt{\det(j_L^*g-\mc{F})},
\end{equation}
where as above $N$ denotes the number of branes in the stack and $\mc{F}=F-j_L^*B$ with $F$ the field strength of the $U(1)$ gauge theory living on the brane. In particular, $d(\mc{F}+j^*_LB)=0$. Moreover, D--branes act as sources for the R-R--fields and couple to these via $\mu\int_L C^{(p+1)}$, where $\mu$ is the brane tension. This term is referred to as the {\em Wess--Zumino action} and in general looks like
\begin{equation}\label{phys_wz}
I_\mathrm{ZW}(L,F)=N\mu\int_L e^{-\mc{F}}\wedge j_L^*C.
\end{equation}

Note that type~IIA and IIB supergravity are related under {\em T--duality}.
Under this transformation, the metric and the B--field of type IIA and IIB are mixed according to the {\em Buscher rules}~\cite{bu87}
(see next section). Moreover, with respect to D--branes, T--duality exchanges the even and odd R-R--potentials and modifies the boundary conditions~(\ref{rumatrix}) in such a way that Dirichlet-- and von~Neumann--conditions are exchanged. This leads to an exchange of D--branes of odd and even dimension.
Whether the dimension of a brane increases or decreases depends on the direction in which the T--duality transformation is carried out. If T--duality is performed along a direction tangent to the brane, the dimension of the transformed brane will decrease, while a transformation along a transverse direction increases the dimension. The mathematical implementation of this will occupy us next.
%
%
%
%
%
\section{\bf T--duality}
\label{tduality}
%
\subsection{The Buscher rules}
\label{buscherrules}

\paragraph{T--dual generalised Riemannian metrics.}
In this section we show how the Buscher rules arise in the generalised context. This parallels work in the generalised complex case~\cite{bb04},~\cite{je04}. 

Consider the vector space $T\oplus T^*\cong\R^n\oplus\R^{n*}$ with a metric splitting $V^+\oplus V^-$ corresponding to $(g,B)$. Further, consider a decomposition $T=\R X\oplus \mc{V}$, where $X$ is a non--zero vector. We define the $1$--form $\theta$ by $\theta(X)=1$, $\ker\theta=\mc{V}$ and obtain the element
\begin{equation}\nonumber
\widetilde{\mc{M}}=X\oplus-\theta\in Pin(n,n).
\end{equation}
Projecting down to $O(n,n)$ gives
\begin{equation}\nonumber
\mc{M}=(Id_{\mc{V}}+\theta^{\dagger})\oplus(Id_{N^*\R X}+X^{\dagger})=\left(\begin{array}{cccc}
\mb{Id_{\mc{V}}}& &\mb{0} &\\ & 0 & & 1\\
\mb{0} & & \,\,\,\mb{Id_{N^*\R X}} & \\
& 1 & & 0
\end{array}\right)=\left(\begin{array}{cc}
\mc{A}& \mc{B} \\ \mc{C} & \mc{D}\end{array}\right).
\end{equation}
Here $\theta^{\dagger}:T\to T^*$ sends $X$ to $\theta$, $\theta^{\dagger}(\mc{V})=0$, and $X^{\dagger}:T^*\to T$ sends $\theta$ to $X$, $X^{\dagger}(N^*\R X)=0$.

\begin{definition}
The {\em T--dual generalised Riemannian metric} is given by $V^{+\top}=\mc{M}(V^+)$.
\end{definition}

\noindent In terms of the corresponding pair $(g^{\top},B^{\top})$, the transformation can be expressed in adapted coordinates as follows.

\begin{prp}\label{buschertrafo}
Extend $X$ to a basis $x_1,\ldots,x_{n-1},x_n=X$ with $x_i\in \mc{V}$, $i=1,\ldots,n-1$, whose dual basis is $x^1,\ldots,x^{n-1},x^n=\theta$. Let $g_{kl}$ and $B_{kl}$ be the coefficients of $g$ and $B$ with respect to this basis. Then the coefficients of $g^{\top}$ and $B^{\top}$ are given according to the Buscher transformation rules (cf. also~\cite{ha00a},~\cite{ha00b} or the appendix in~\cite{kstt03}), namely
\begin{equation}\nonumber
\begin{array}{ll}
g_{kl}^{\top}=g_{kl}-\frac{1}{g_{nn}}(g_{kn}g_{nl}-B_{kn}B_{nl}),\quad & g^{\top}_{kn}=\frac{1}{g_{nn}}B_{kn},\quad g^{\top}_{nn}=\frac{1}{g_{nn}}\\[5pt]
B_{kl}^{\top}=B_{kl}+\frac{1}{g_{nn}}(g_{kn}B_{ln}-B_{kn}g_{ln}),\quad & B^{\top}_{kn}=\frac{1}{g_{nn}}g_{kn}.
\end{array}
\end{equation}
\end{prp}

\begin{prf}
The elements of $V^{+\top}$ are given by $Y\oplus P^{+\top}Y=\mc{A}X+\mc{B}P^+X\oplus\mc{C}X+\mc{D}P^+X$ for $X\oplus P^+X\in V^+$, hence $P^{+\top}=(\mc{C}+\mc{D}P^+)(\mc{A}+\mc{B}P^+)^{-1}$. Writing
\begin{equation}\nonumber
P^+=g+B=\left(\begin{array}{cc}\underline{g} & \alpha\\\alpha^{tr} & q\end{array}\right)+\left(\begin{array}{cc}\underline{B} & \beta\\-\beta^{tr} & 0\end{array}\right),
\end{equation}
we obtain
\begin{equation}\nonumber
(\mc{A}+\mc{B}P^+)^{-1}=\left(\begin{array}{cc}\underline{Id} & 0\\ -(\alpha-\beta)^{tr}/q & 1/q\end{array}\right)
\end{equation}
and
\begin{equation}\nonumber
(\mc{C}+\mc{D}P^+)=\left(\begin{array}{cc}\underline{g}+\underline{B} & \alpha+\beta\\0 & 1\end{array}\right).
\end{equation}
Consequently,
\begin{equation}\nonumber
P^{+\top}=\left(\begin{array}{cc}\underline{g}+\underline{B}-(\alpha+\beta)(\alpha-\beta)^{tr}/q & (\alpha+\beta)/q\\-(\alpha-\beta)^{tr}/q & 1/q\end{array}\right)=g^{+\top}+B^{+\top}, 
\end{equation}
that is,
\begin{equation}\nonumber
g^{+\top}=\left(\begin{array}{cc}\underline{g}-(\alpha\alpha^{tr}-\beta\beta^{tr})/q & \beta/q\\\beta^{tr}/q & 1/q\end{array}\right)
\end{equation}
and
\begin{equation}\nonumber B^{+\top}=\left(\begin{array}{cc}\underline{B}+(\alpha\beta^{tr}-\beta\alpha^{tr})/q & \alpha/q\\-\alpha^{tr}/q & 0\end{array}\right),
\end{equation}
whence the assertion.
\end{prf}

\paragraph{T--dual spinors.}
For $\rho\in S_{\pm}$, its {\em T--dual spinor} is
\begin{equation}\nonumber
\rho^{\top}=\widetilde{\mc{M}}\bullet\rho.
\end{equation}
If the collection of chiral spinors $\rho_1,\ldots,\rho_s$ defines a generalised $G_L\times G_R$--structure, that is, their stabiliser in $Spin(n,n)_+$ is $G_L\times G_R$, the collection $\rho_1^{\top},\ldots,\rho_s^{\top}$ induces the {\em T--dual} $G_L\times G_R$--structure associated with the conjugated group $\widetilde{\mc{M}}(G_L\times G_R)\widetilde{\mc{M}}$. In particular, equivariance implies the induced generalised Riemannian metric to be given by $(g^{\top},B^{\top})$, whence $\gtilde^{\top}\rho^{\top}=(\gtilde\rho)^{\top}$.

\begin{rmk}
Note that T--duality reverses the parity of the spinor. This reflects the physical fact that T--duality exchanges type IIA with type IIB theory and hence maps odd R-R potentials into even R-R potentials (cf. Section~\ref{physics}).
\end{rmk}

\noindent The operator $\widetilde{\mc{M}}$ also maps pure spinors to pure spinors. More concretely, if $\tau$ is pure with annihilator $W_{\tau}$, then $\tau^{\top}=\widetilde{\mc{M}}\bullet\tau$ is pure with annihilator $\mc{M}(W_{\tau})$. Hence, T--duality associates to any isotropic pair $(L,F)$ a uniquely determined T--dual isotropic pair $(L^{\top},F^{\top})$ specified by the condition $[\tau^{\top}]=[\tau_{L^{\top},F^{\top}}]$. Further, $\widetilde{\mc{M}}$ is norm--preserving in the sense that for $\tau\in S_{\pm}$, we have $\eunorm{\tau^{\top}}_{\gtilde^{\top}}=\eunorm{\tau}_{\gtilde}$. Therefore
\begin{equation}\label{tdualcalcond}
\langle\rho,\tau\rangle\leq\eunorm{\tau}_{\gtilde}\quad\mbox{if and only if}\quad\langle\rho^{\top},\tau^{\top}\rangle\leq\eunorm{\tau^{\top}}_{\gtilde^{\top}},
\end{equation}
and equality holds on the left hand side if and only if equality holds on the right hand side.

\begin{prp}\label{tdualrank}
A spinor $\rho\in S_{\pm}$ defines a calibration for the generalised Riemannian metric $(g,B)$ on $T\oplus T^*$ if and only if $\rho^{\top}=\widetilde{\mc{M}}\bullet\rho$ defines a calibration for the T--dual generalised metric $(g^{\top},B^{\top})$. A (calibrated) isotropic pair $(L,F)$ corresponds to a uniquely (calibrated) isotropic pair $(L^{\top},F^{\top})$. The rank of $\tau_{L^{\top},F^{\top}}$ is given as follows:
\begin{itemize}
	\item If $X\in L$, then $\rk\tau_{L^{\top},F^{\top}}=\rk\tau_{L,F}-1$.
	\item If $X\not\in L$, then $\rk\tau_{L^{\top},F^{\top}}=\rk\tau_{L,F}+1$.
\end{itemize}
\end{prp}

\begin{prf}
In virtue of the preceding discussion, only the last assertion requires proof. By definition of the rank (cf. Proposition~\ref{rankdef}),
\begin{eqnarray*}
\rk\tau_{L^{\top},F^{\top}} & = & n-\dim\big(\mc{M}(W_{L,F})\cap T^*\big)\\
& = & n-\dim\mc{M}\big(W_{L,F}\cap\mc{M}(T^*)\big)\\
& = & n-\dim\big(W_{L,F}\cap (\R X\oplus N^*\R X)\big)\\
& = & n-\dim\big(W_L\cap (\R X\oplus N^*\R X)\big),
\end{eqnarray*}
where for the last line we used that $X\llcorner F\in N^*\!L$ if and only if $X\llcorner F=0$. 

Now suppose that $X\in L$. In this case, $L=\R X\oplus L\cap\mc{V}$, hence we may choose a basis $x_1,\ldots,x_n$ of $T$ such that $x_1=X$, $x_2,\ldots,x_k$ generates $L\cap\mc{V}$, and $x_2,\ldots,x_n$ generates $\mc{V}$. Let $x^1,\ldots,x^n$ denote the dual basis. Hence $\theta=x^1$ and
\begin{equation}\nonumber
W_L=L\oplus N^*\!L=\{a^1X+\sum_{i=2}^k a^ix_i\oplus\sum_{j=k+1}^nb_ix^i\}.
\end{equation}
The intersection with $\R X\oplus N^*\R X$ is spanned by $X,x^{k+1},\ldots,x^n$, so that $\rk\tau_{L^{\top},F^{\top}}=n-k-1=\rk\tau_{L,F}-1$. 

Next assume that $X\not\in L$. Then there exists a basis $x_1=X,x_2,\ldots,x_n$ of $\R X\oplus \mc{V}$ such that $L$ is spanned by $cX+x_2,\ldots,x_{k+1}$ (possibly $c=0$). Again $\theta=x^1$, and
\begin{equation}\nonumber
W_L=L\oplus N^*\!L=\{a^1(cX+x_2)+\sum_{i=3}^{k+1} a^ix_i\oplus b(\theta-cx^2)+\sum_{j=k+2}^nb_ix^i\}.
\end{equation}
The intersection with $\R X\oplus N^*\R X$ is now spanned by $x^{k+2},\ldots,x^n$, hence  $\rk\tau_{L^{\top},F^{\top}}=n-k+1=\rk\tau_{L,F}+1$.
\end{prf}

\begin{ex}
Consider a $G_2$--structure on $T^7$. As we learned from the examples in Section~\ref{calsub}, an isotropic pair $(L,F)$ consisting of a coassociative subplane $L$ and an anti--self--dual $2$--form $F$ is calibrated by $\rho=1-\star\varphi$. Now choose a vector $X\in L$ and a complement $\mc{V}$ so that $T=\R X\oplus\mc{V}$. Then $(L^{\top},F^{\top})$ is a calibrated pair with $L^{\top}$ of dimension $\dim L-1=3$. On the other hand, choosing a vector $X\not\in L$ yields a T--dual calibrated pair with a $5$--dimensional plane $L^{\top}$. In particular, we see that calibrated isotropic pairs of ``exotic'' dimension (in comparison to the straight case, e.g. $3$ and $4$ for $G_2$) can occur.
\end{ex}
%
\subsection{Integrability}
In this section, we discuss a local version of T--duality over $M=\R^n$ (now seen as a manifold) endowed with a generalised metric $(g,H=dB)$. The induced generalised tangent bundle $\mb{E}(H)\to\R^n$ is equivalent as a vector bundle to $T\R^n\oplus T^*\R^n$, but inequivalent from a generalised point of view if $H\not=0$, as this implies twisting with a {\em non--closed} $B$--field (cf. also Section~\ref{hfluxtwist}). Fix a dilaton field $\phi$ so that $\nu=\exp(2\phi)\nu_g$. Further, let $X$ be a nowhere vanishing vector field transversal to the $n-1$--dimensional distribution $\mc{V}$ so that $T_p\R^n=\R X(p)\oplus\mc{V}(p)$. This determines the $1$--form $\theta$ and hence the operator $\widetilde{\mc{M}}=X\oplus-\theta$ on $\mb{S}(\mb{E})$. T--duality induces a map $\widetilde{\mc{M}}:\Omega^{ev,od}(\R^n)\to\Omega^{od,ev}(\R^n)\!\cong\!\mb{S}_{\mp}\big(\mb{E}(H^{\top})\big)$ (where the last isomorphism is induced by $\mf{L}^{\nu\top}$, the isomorphism between spinors and differential forms associated with $H^{\top}=dB^{\top}$) defined by
\begin{equation}\nonumber
\widetilde{\mc{M}}(\alpha)=e^{B^{\top}}\!\wedge\widetilde{\mc{M}}\bullet e^B\bullet\rho.
\end{equation}
The {\em T--dual dilaton field} $\phi^{\top}$ is defined by 
\begin{equation}\nonumber
\nu=\exp(2\phi^{\top})\nu_{g^{\top}}.
\end{equation}

\begin{prp}
The dilaton transforms under the Buscher rule
\begin{equation}\nonumber
\phi^{\top}=\phi-\ln\eunorm{X}.
\end{equation}
\end{prp}

\begin{prf}
Let $q=\eunorm{X}^2$. It suffices to show that $\nu_g=q\cdot\nu_{g^{\top}}$. This is a tensorial identity and can be computed pointwise. At a given point $p$, we may assume that $\mc{V}$ is spanned by $\partial_{x^i}(p)$, $i\leq n-1$, and $X(p)=\partial_{x^n}$ for coordinates $x^1,\ldots,x^n$. By Proposition~\ref{buschertrafo}, $g^{\top}(p)=A^{tr}\circ g(p)\circ A$, where $A\in GL(n)$ is given by
\begin{equation}\nonumber
A=\left(\begin{array}{cc}\underline{Id}& 0\\(\beta-\alpha)^{tr}/q & 1/q\end{array}\right).
\end{equation}
The claim follows from $\nu_{g^{\top}}(p)=\det A^{-1}\cdot\nu_g(p)$. 
\end{prf}

\begin{ex}
Consider the spinor field $\rho$ defined by $\mf{L}^{\nu_g}(\rho)=[\Psi\otimes\Psi]^{od}$, where $\Psi$ comes from an ordinary $G_2$--structure. The T--dual spinor field $\rho^{\top}$ induces also a generalised $G_2$--structure and is therefore determined by $\mf{L}^{\nu_g\,\top}(\rho^{\top})=\mf{L}^{\phi^\top}(\rho^{\top})=e^{-\phi^{\top}}[\Psi^{\top}_L\otimes\Psi^{\top}_R]^{od}$.
\end{ex}

Under additional assumptions, T--duality also preserves closeness of spinor fields.

\begin{thm}\label{tdualint}
Assume that $(g,B,\phi)$ are $X$--invariant, that is $\mc{L}_{\!X}g=0$, $\mc{L}_X\!B=0$ and $\mc{L}_X\phi=0$. Further, assume that $\mc{L}_X\theta=0$. 
If $\rho$ and $\varphi$ are two $\mb{E}(H)$--spinor fields such that $\mc{L}_X\mf{L}^{\nu}(\rho)$ and $\mc{L}_X\mf{L}^{\nu}(\varphi)=0$, then $d_{\nu}\rho=\varphi$ if and only if $d_{\nu}\rho^{\top}=-\varphi^{\top}$, that is
\begin{equation}\label{firstint}
d_H\mf{L}^{\nu}(\rho)=\mf{L}^{\nu}(\varphi)\quad\mbox{if and only if}\quad d_{H^{\top}}\mf{L}^{\nu\,\top}(\rho^{\top})=-\mf{L}^{\nu\,\top}(\varphi^{\top}).
\end{equation}
\end{thm}

\begin{prf}
Fix coordinates $x^1,\ldots,x^n$ so that $\nu=e^{2\widetilde{\phi}_{\underline{x}}}\partial_{x^1}\wedge\ldots\wedge\partial_{x^n}$ with $\widetilde{\phi}_{\underline{x}}=\phi-\ln(\det_{\underline{x}}g)/4$. The subscript $\underline{x}$ reminds us that this quantity depends upon the choice of coordinates.  It is invariant under T--duality in the sense that $\widetilde{\phi}_ {\underline{x}}=\phi^{\top}-\ln(\det_x g^{\top})/4=\widetilde{\phi}_ {\underline{x}}^{\top}$. 
By definition, $d_{\nu}\rho=\varphi$ is then equivalent to
\begin{equation}\label{tint1}
de^{-\widetilde{\phi}_{\underline{x}}}\rho=e^{-\widetilde{\phi}_{\underline{x}}}\varphi.
\end{equation}
Let 
\begin{equation}\nonumber
\rho=\rho_0+\theta\wedge\rho_1,\quad\varphi=\varphi_0+\theta\wedge\varphi_1
\end{equation}
be the decomposition into basic forms, that is, $X\llcorner\rho_{0,1}=0$ and $X\llcorner\varphi_{0,1}=0$.
The assumptions imply that $\mc{L}_X\mf{L}^{\nu}(\rho)=e^{-\phi}\exp(-B)\wedge\mc{L}_X(\rho)=0$, whence $X\llcorner d\rho_{0,1}=0$. Analogously, $X\llcorner d\varphi_{0,1}=0$ holds. Therefore,~(\ref{tint1}) is equivalent to
\begin{equation}\label{tint}
\varphi_{0}=-d\widetilde{\phi}_ {\underline{x}}\wedge\rho_{0}+d\rho_{0}+d\theta\wedge\rho_1,\quad\varphi_{1}=d\widetilde{\phi}_ {\underline{x}}\wedge\rho_{1}-d\rho_{1}.
\end{equation}
For the right hand side of~(\ref{firstint}), we find
\begin{equation}\label{tint2}
-\varphi^{\top}_0=-d\widetilde{\phi}_ {\underline{x}}^{\top}\wedge\rho_0^{\top}+d\rho^{\top}_0+d\theta\wedge\rho_1^{\top},\quad-\varphi_1^{\top}=d\widetilde{\phi}_ {\underline{x}}^{\top}\wedge\rho^{\top}_1-d\rho^{\top}_1,
\end{equation}
where
\begin{equation}\nonumber
\rho^{\top}=\rho_0^{\top}+\theta\wedge\rho_1^{\top}=\widetilde{\mc{M}}\bullet\rho=-(X\llcorner+\theta\wedge)(\rho_{0}+\theta\wedge\rho_{1})=-\rho_{1}-\theta\wedge\rho_{0},
\end{equation}
and similarly for $\varphi^{\top}$. As $\widetilde{\phi}_ {\underline{x}}^{\top}=\widetilde{\phi}_ {\underline{x}}$, we deduce from this that~(\ref{tint2}) is equivalent to~(\ref{tint}) and thus to~(\ref{tint1}).
\end{prf}

\begin{ex}
To illustrate the theorem we take up the previous example where in addition we assume: (i) the $G_2$--invariant $3$--form $\varphi$ is closed and coclosed (ii) there exists a vector field $X$ of non--constant length such that $\mc{L}_X\varphi=0$. In particular, $X$ defines a Killing vector field. Local examples of this exist in abundance, cf. for instance~\cite{apsa04}. Since $\mf{L}^{\nu}(\rho)=[\Psi\otimes\Psi]^{od}$ is closed, so is the T--dual $\mf{L}^{\nu\,\top}(\rho^{\top})=e^{-\phi^{\top}}[\Psi^\top_L\otimes\Psi^\top_R]^{od}$. Since $\varphi$ is also co--closed, the same holds for the even spinor field $\gtilde(\rho)$ and its T--dual. Hence
\begin{equation}\label{susyg2}
d_{H^{\top}}e^{-\phi^{\top}}[\Psi^\top_L\otimes\Psi^\top_R]=0,
\end{equation} 
that is, the generalised $G_2$--structure is {\em integrable}~\cite{jewi05a},~\cite{wi04}, where equation~(\ref{susyg2}) was shown to encode the supersymmetry variations of type II supergravity. In particular, this implies $H^{\top}\not=0$, for $\phi^{\top}=-ln\eunorm{X}$ is not constant. Although we started with $B=0$ and $\phi=0$, the Buscher rules imply that we acquire a generalised $G_2$--structure with non--trivial $B$-- and dilaton field on the T--dual side.
\end{ex}

Since being a calibration is a pointwise condition, it follows from~(\ref{tdualcalcond}) that this property is preserved under T--duality. Similarly, one would expect calibrated isotropic pairs to transform under T--duality in the vein of Proposition~\ref{tdualrank}. However, two problems may occur: Firstly, if $(L,\mc{F})$ is calibrated by $\rho$, the rank of $\tau^{\top}_{L,\mc{F}}$ is not constant in general, which prevents this spinor field to be induced by an isotropic pair. In view of Proposition~\ref{tdualrank}, we either need to restrict $L$ to the open subset for which $X\not\in T_pL$, or to assume $X\in\Gamma(TL)$. Secondly, if there is a pair $(L',\mc{F}')$ such that $[\tau^{\top}_{L,\mc{F}}]=[\tau_{L',\mc{F}'}]$, it might not be isotropic, that is, we possibly have $d\mc{F}'+j^*_{L'}H^{\top}\not=0$. This requires to rephrase the isotropy condition on the pair in terms of integrability conditions of the induced pure spinor field.

\begin{definition}
Let $\mb{E}\to M$ be a generalised tangent bundle and $\nu$ a trivialisation of $\Lambda^nTM$. A pure spinor field half--line in $\mb{S}(\mb{E})_{\pm}$ is called {\em integrable} if and only if it admits a $d_{\nu}$--closed representative of constant rank.
\end{definition}

\noindent By abuse of language, we also refer to the representative as integrable.

\begin{lem}
Let $(L,\mc{F})$ be an isotropic pair for $(\R^n,H=dB)$. Then locally, $[\tau_{L,\mc{F}}]$ is induced by an integrable half--line, that is there exists an open set $U\subset\R^N$ and an integrable pure spinor field $\tau\in\mb{S}(\mb{E})_{|U}$ with $[\tau_{|U\cap L}]=[\tau_{L,\mc{F}|U\cap L}]$. Furthermore, if $\mc{L}_XB=0$ and $X$ is a vector field either transversal to $L$, or contained in $TL$ with $\mc{L}_X\mc{F}=0$, then $\tau$ can be constructed to satisfy $\mc{L}_X\mf{L}^{\nu}(\tau)=0$.

Conversely, an integrable pure spinor field $\tau$ over $\R^n$ gives locally rise to a foliation into isotropic pairs $(L_{\underline{r}},\mc{F}_{\underline{r}})$ such that $[\tau_{|L_{\underline{r}}}]=[\tau_{L_{\underline{r}},\mc{F}_{\underline{r}}}]$
\end{lem}

\begin{prf}
Let us start with the converse. Restricting $\tau$ to some open subset $U$, we can write $\mf{L}^{\nu}(\tau)=\exp(\mc{F}_{\tau})\wedge\varpi$ with $\mc{F}_{\tau}\in\Omega^2(U)$ and $\omega=\theta^1\wedge\ldots\wedge\theta^{n-k}$, $\theta^i\in T^*U$. Closeness is equivalent to \begin{equation}\label{closedspinor}
d\varpi=0,\quad (d\mc{F}_{\tau}+H)\wedge\varpi=0.
\end{equation}
By Frobenius' theorem, the smooth distribution ${\rm Ann}\,\varpi=\{Y\in\mf{X}(U)\,|\,Y\llcorner\varpi=0\}$ is thus integrable. Hence, on $U$ there exist coordinates $y^1,\ldots,y^n$ (possibly upon shrinking $U$) such that for $\underline{r}\in\R^{n-k}$ and $p\in L_{\underline{r}}=\{p\in U\,|\,y^i(p)=r^i,\,i=k+1,\ldots,n\}$, one has $T_pL_{\underline{r}}={\rm Ann}\,\varpi(p)$. Further by~(\ref{closedspinor}), $d\mc{F}_{\tau}+H\in\mc{N}_{L_{\underline{r}}}^3$, the space of $3$--forms of the ideal in $\Omega^*(\R^n)_{|L_{\underline{r}}}$ generated by $N^*\!L_{\underline{r}}\subset T^*\R^n_{|L_{\underline{r}}}$. Consequently, defining $\mc{F}_{\underline{r}}=j^*_{L_{\underline{r}}}\mc{F}_{\tau}$,
\begin{equation}\nonumber
d\mc{F}_{\underline{r}}+j_{L_{\underline{r}}}^*H=j_{L_{\underline{r}}}^*(d\mc{F}_{\tau}+H)=0.
\end{equation}
By construction, $(\mc{F}_{_{\tau}|L_{\underline{r}}}-\mc{F}_{\underline{r}})\wedge\varpi=0$, hence $[\tau_{|L_{\underline{r}}}]=[\tau_{L_{\underline{r}},\mc{F}_{\underline{r}}}]$ for suitably oriented $L_{\underline{r}}$.

Now consider the isotropic pair $(L,\mc{F})$. We can fix around $p\in L$ a coordinate system $y^1,\ldots,y^n$ on $U\subset\R^n$ with $L=\{y^{k+1}=\ldots=y^n=0\}$. Define the closed $2$--form $F_L=(\mc{F}+j^*_LB)_{|U\cap L}$, which we extend trivially to a $2$--form $F$ on $U$, that is, $F(p)=F(p^1,\ldots,p^n)=F_L(p^1,\ldots,p^k,0,\ldots,0)$. Set $\mc{F}=F-B_{|U}$ and let $\tau\in\mb{S}(\mb{E})$ be determined by $\mf{L}^{\nu}(\tau)=e^{\mc{F}}\wedge dy^{k+1}\wedge\ldots\wedge dy^n$. Then $[\tau_{|U\cap L}]=[\tau_{L,\mc{F}|U\cap L}]$. We claim $\tau$ to be $d_{\nu}$--closed. Indeed, this is equivalent to
\begin{equation}\label{closed}
(d\mc{F}+H_{|U})\wedge dx^{k+1}\wedge\ldots\wedge dx^n=0,
\end{equation} 
hence to $dF\wedge dy^{k+1}\wedge\ldots\wedge dy^n=0$. But this holds on $L$, hence on $U$ for we extended $F$ trivially. In particular, $j^*_{L_{\underline{r}}}(d\mc{F}+H_{|U})=0$ using the previous notation, so the pairs $(L_{\underline{r}},\mc{F}_{L_{\underline{r}}})$ are isotropic. If $X$ is a vector field transversal to $L$, we carry out this construction with a coordinate system which in addition satisfies $\partial_{y^n}=X$. Then $\mf{L}^{\nu}(\tau)$ does not depend on $y^n$ for $\mc{L}_X\mc{F}=0$, hence $\mc{L}_X\mf{L}^{\nu}(\tau)=0$. On the other hand, if $X$ restricts to a vector field on $L$, we can fix a coordinate system with $\partial_{y^1}=X$ and $L=\{y^{k+1}=\ldots=y^n=0\}$. Again, $\tau$ is $X$--invariant as $\mc{L}_X\mc{F}=0$.
\end{prf}

\noindent We are now in a position to prove

\begin{thm}\label{caltdual}\hfill\newline
{\rm\noindent(i)}  The spinor $\rho\in \Gamma\big(\mb{S}(\mb{E}(H))_{\pm}\big)$ defines a calibration for $(g,H=dB)$ if and only if $\rho^{\top}\in\Gamma\big(\mb{S}(\mb{E}(H^{\top}))_{\pm}\big)$ does so for $(g^{\top},H^{\top}=dB^{\top})$. Furthermore, if $d_{\phi}\rho=d_{\phi}\gamma$ under the assumptions of Theorem~\ref{tdualint}, then the calibrated isotropic pairs for $\rho$ and $\rho^{\top}$ extremalise the functionals $\mc{E}_{\phi,\gamma}$ and $\mc{E}_{\phi^{\top},-\gamma^{\top}}$ respectively. 

{\rm\noindent(ii)} Under the assumptions of (i), let $(L,\mc{F})$ be a calibrated isotropic pair. If $X(p)\not\in T_pL$ for some $p\in L$, there exists a calibrated isotropic pair $(L^{\top},\mc{F}^{\top})$ with
\begin{itemize}
	\item $p\in L^{\top}$ and $\dim L^{\top}=\dim L+1$
	\item $[\tau_{L^{\top},\mc{F}^{\top}}(p)]=[\tau^{\top}_{L,\mc{F}}(p)]$.
\end{itemize}
On the other hand, if $X\in\Gamma(TL)$ and $\mc{L}_X\mc{F}=0$, there exists a calibrated isotropic pair $(L^{\top},\mc{F}^{\top})$ with
\begin{itemize}
	\item $p\in L^{\top}$ and $\dim L^{\top}=\dim L-1$
	\item $[\tau_{L^{\top},\mc{F}^{\top}}(p)]=[\tau^{\top}_{L,\mc{F}}(p)]$.
\end{itemize} 
\end{thm}

\begin{rmk}
The dimension shift coincides with expectations from physics, cf. Section~\ref{physics2}.
\end{rmk}

\begin{prf}
In virtue of the preceding discussion and Theorem~\ref{tdualint}, only assertion (ii) requires proof. By the previous lemma, we can locally extend $[\tau_{L,\mc{F}}]$ to an integrable pure spinor field half--line $[\tau]$ with $\mc{L}_X\mf{L}^{\phi}(\tau)=0$ and $d_{\phi}\tau=0$ (on $U$, the domain of $\tau$). Consequently, $d_{\phi^{\top}}\tau^{\top}=0$ by Theorem~\ref{tdualint}. If $X$ is transversal to $L$, then $X\not\in TL_{\underline{r}}$. Hence $\tau^{\top}$ is of constant rank $\rk(\tau)+1$ and therefore integrable on $U$. Appealing to the converse of the lemma, we deduce the existence of the isotropic pair $(L^{\top},\mc{F}^{\top})$. The rest follows analogously.
\end{prf}

\begin{ex}
We saw that if $(M,\varphi)$ is a classical $G_2$--manifold with $G_2$--invariant spinor $\Psi$, all calibrated submanifolds $L$ give rise to calibrated isotropic pairs $(L,\mc{F}=0)$. If $\mc{L}_X\varphi=0$  and $X$ is transversal or restricts to a vector field on $L$, we obtain a calibrated isotropic pair $(L^{\top},\mc{F}^{\top})$ for the non--straight generalised $G_2$--structure $[\Psi\otimes\Psi]^{od\,\top}$ by choosing some distribution complementary to $X$ with $\mc{L}_X\theta=0$, which locally is always possible.
\end{ex}

\bigskip

\noindent
F.~\textsc{Gmeiner}: NIKHEF, Kruislaan 409, 1098 DB Amsterdam, The Netherlands\\ 
e--mail: \texttt{fgmeiner@nikhef.nl}

\medskip

\noindent
F.~\textsc{Witt}: NWF I -- Mathematik, Universit\"at Regensburg,
D--93040 Regensburg, F.R.G.\\
e--mail: \texttt{Frederik.Witt@mathematik.uni-regensburg.de}

\addcontentsline{toc}{section}{References}
\begin{thebibliography}{10}
\bibitem{apsa04}
V.~Apostolov and S.~Salamon:
{\em K{\"a}hler reduction of metrics with holonomy $G_2$.} 
Comm. Math. Phys. {\bf 246}(1):43--61 (2004).


\bibitem{bb04}
O.~Ben--Bassat:
{\em Mirror symmetry and generalized complex manifolds.} 
J. Geom. Phys. {\bf 56}:533--558 (2006).


\bibitem{bebo04}
O.~Ben--Bassat and M.~Boyarchenko:
{\em Submanifolds of generalized complex manifolds.} 
J. Symplectic Geom. {\bf 2}(3):309--355 (2004).


\bibitem{begr06}
I.~Benmachiche and T.~Grimm: 
{\em Generalized $\mc{N}=1$ orientifold compactifications and the Hitchin functionals.} Nucl. Phys. B {\bf 748}:200--252 (2006).


\bibitem{bkop97}
E.~Bergshoeff, R.~Kallosh, T.~Ortin and G.~Papadopoulos:
{\em kappa--symmetry, supersymmetry and intersecting branes.} 
Nucl. Phys. B {\bf 502}:149--169 (1997).


\bibitem{bkop01}
E.~Bergshoeff, R.~Kallosh, T.~Ortin, D.~Roest and A.~Van~Proeyen:
{\em New formulations of $D=10$ supersymmetry and $D8-O8$ domain walls.}
Class. Quant. Grav. {\bf 18}:3359--3382 (2001).


\bibitem{bem04}
P.~Bouwknegt, J.~Evslin and V.~Mathai:
{\em T--Duality: Topology Change from $H$--flux.}
Comm. Math. Phys. {\bf 249}(2):383--415 (2004).


\bibitem{brs05}
U.~Bunke, P.~Rumpf and T.~Schick:
{\em The topology of T--duality for $T^n$--bundles}, math/0501487.


\bibitem{busc04}
U.~Bunke and T.~Schick:
{\em On the topology of T--duality.} 
Rev. Math. Phys. {\bf 17}:77--112 (2005).


\bibitem{bu87}
T.~Buscher:
{\em A symmetry of the string background field equations.} Phys. Lett. B {\bf 194}:59--62 (1987).


\bibitem{ch96}
C.~Chevalley:
{\em The algebraic theory of spinors and Clifford algebras.}
Collected works vol. {\bf 2}, Springer, Berlin (1996).


\bibitem{cgj04}
S.~Chiantese, F.~Gmeiner and C.~Jeschek:
{\em Mirror symmetry for topological sigma models with generalized K\"ahler geometry.} Int. J. Mod. Phys. A {\bf 21}:2377--2390 (2006).


\bibitem{co90}
T.~Courant:
{\em Dirac manifolds.}
Trans. Amer. Math. Soc. {\bf 319}:631--661 (1990).


\bibitem{daha93}
J.~Dadok and R.~Harvey:
{\em Calibrations and spinors.}
Acta Math. {\bf 170}(1):83--120 (1993).


\bibitem{gmw04}
J.~Gauntlett, D.~Martelli and D.~Waldram: 
{\em Superstrings with intrinsic torsion.}
Phys.~Rev.~D {\bf 69}:086002 (2004).


\bibitem{gip02}
J.~Gutowski, S.~Ivanov and G.~Papadopoulos:
{\em Deformations of generalized calibrations and compact non--Kahler manifolds with vanishing first chern class.} 
Asian J. Math. {\bf 7}(1):39--79 (2003).


\bibitem{gupa99}
J.~Gutowski and G.~Papadopoulos:
{\em AdS calibrations.} 
Phys. Lett. B {\bf 462}:81--88 (1999).


\bibitem{ha91}
R.~Harvey:
{\em Spinors and Calibrations.}
Perspectives in Mathematics vol. {\bf 9}, Academic Press, Boston (1990).


\bibitem{hala82}
R.~Harvey and H.~Lawson:
{\em Calibrated geometries.} 
Acta Math. {\bf 148}:47--157 (1982).


\bibitem{ha00a}
S.~Hassan:
{\em T--duality, space--time spinors and R-R fields in curved backgrounds.}
Nucl. Phys. B {\bf 568}:145--161 (2000).


\bibitem{ha00b}
S.~Hassan:
{\em $SO(d,d)$ transformations of Ramond--Ramond fields and space--time spinors.}
Nucl. Phys. B {\bf 583}:431--453 (2000).


\bibitem{hi03}
N.~Hitchin:
{\em Generalized Calabi--Yau manifolds.}
Quart. J. Math. Oxford Ser. {\bf 54}:281--308 (2003).


\bibitem{hi05}
N.~Hitchin:
{\em Brackets, forms and invariant functionals.} Asian J. Math. {\bf10}(3):541--560 (2006).


\bibitem{je04}
C.~Jeschek:
{\em Generalized Calabi--Yau structures and mirror symmetry,} hep-th/0406046.


\bibitem{jewi05a}
C.~Jeschek and F.~Witt:
{\em Generalised ${G}_2$--structures and type IIB superstrings.}
JHEP {\bf 0503}:053 (2005).


\bibitem{jewi05b}
C.~Jeschek and F.~Witt:
{\em Generalised geometries, constrained critical points and Ramond--Ramond fields}, math.DG/0510131.


\bibitem{jo03}
C.~Johnson:
{\em D--branes.}
Cambridge University Press, Cambridge (2003).


\bibitem{jo04}
J.~Joyce:
{\em The exceptional holonomy groups and calibrated geometry,} math.DG/0406011.


\bibitem{kstt03}
S.~Kachru, M.~Schulz, P.~Tripathy and S.~Trivedi:
{\em New supersymmetric string compactifications.}
JHEP {\bf 0303}:061 (2003).


\bibitem{kali04}
A.~Kapustin and Y.~Li:
{\em Topological sigma--models with H--flux and twisted generalized complex manifolds.} Adv. Theor. Math. Phys. {\bf11}(2):261--290 (2007).

  
\bibitem{ko05}
P.~Koerber:
{\em Stable D--branes, calibrations and generalized Calabi--Yau geometry.}
JHEP {\bf 0508}:099 (2005).


\bibitem{mmms00}
M.~Marino, R.~Minasian, G.~Moore and A.~Strominger:
{\em Nonlinear instantons from supersymmetric $p$--branes.}
JHEP {\bf 0001}:005 (2000).


\bibitem{ma06}
L.~Martucci:
{\em D--branes on general $\mc{N}=1$ backgrounds: Superpotentials and D--terms.}
JHEP {\bf 0511}:048 (2005).


\bibitem{masm05}
L.~Martucci and P.~Smyth:
{\em Supersymmetric D--branes and calibrations on general $\mc{N}=1$ backgrounds.}
JHEP {\bf 0511}:048 (2005).

    
\bibitem{po96}
J.~Polchinski:
{\em Lectures on D--branes,} hep-th/9611050.


\bibitem{syz96}
A.~Strominger, S.-T.~Yau and E.~Zaslow:
{\em Mirror symmetry is T--duality.}
Nucl. Phys. B {\bf 479}:243--259 (1996).


\bibitem{wa89}
M.~Wang:
{\em Parallel spinors and parallel forms.}
Ann. Global Anal. Geom. {\bf 7}(1):59--68 (1989).


\bibitem{wi04}
F.~Witt:
{\em Generalised $G_2$--manifolds.}
Comm. Math. Phys. {\bf 265}(2):275--303 (2006).


\bibitem{wi05}
F.~Witt:
{\em Special metric structures and closed forms.}
DPhil thesis University of Oxford (2005).


\bibitem{wi06}
F.~Witt:
{\em Metric bundles of split signature and type II supergravity.}
To appear in: ``Recent developments in pseudo-Riemannian Geometry'' ed. by D.~Alekseevski and H.~Baum, ESI--Series on Mathematics and Physics.


\bibitem{za04}
M.~Zabzine:
{\em Geometry of D--branes for general $\mc{N}=(2,2)$ sigma models.}
Lett. Math. Phys. {\bf 70}:211--221 (2004).
\end{thebibliography}
\end{document}